\documentclass[10pt,reqno]{amsart}
\usepackage{amssymb,amsthm,amsmath,latexsym,fullpage}
\usepackage{amsfonts,subfigure,graphicx}
\newcommand{\N}{{I\!\! N}}
\newcommand{\R}{{I\!\! R}}
\newcommand{\cal}{\mathcal}

\numberwithin{equation}{section}

\begin{document}
\title{Existence and Stability of Compressible Current-Vortex Sheets
 in Three-Dimensional Magnetohydrodynamics}

\author{Gui-Qiang Chen}
\address{G.-Q. Chen,
Institute of Mathematics and Key Laboratory of Mathematics for
Nonlinear Science, Fudan University, Shanghai 200433, PRC;
Department of Mathematics, Northwestern University, 2033 Sheridan
Road, Evanston, Illinois 60208, USA}
\email{gqchen@math.northwestern.edu}
\urladdr{http://www.math.northwestern.edu/$\sim$gqchen}
\author{Ya-Guang Wang}
\address{Y.-G. Wang, Department of Mathematics,
     Shanghai Jiao Tong University, 200240 Shanghai, PRC}
\email{E-mail: ygwang@sjtu.edu.cn}

\date{September 8, 2006}

\begin{abstract}
Compressible vortex sheets are fundamental waves, along with shock
and rarefaction waves, in entropy solutions to the multidimensional
hyperbolic systems of conservation laws; and understanding the
behavior of compressible vortex sheets is an important step towards
our full understanding of fluid motions and the behavior of entropy
solutions.
For the Euler equations in two-dimensional gas dynamics, the
classical linearized stability analysis on compressible vortex
sheets predicts stability when the Mach number $M>\sqrt{2}$ and
instability when $M<\sqrt{2}$; and Artola-Majda's analysis reveals
that the nonlinear instability may occur if planar vortex sheets are
perturbed by highly oscillatory waves even when $M>\sqrt{2}$.
For the Euler equations in three-dimensions,
every compressible vortex sheet is violently unstable and this
violent instability is the analogue of the Kelvin-Helmholtz
instability for incompressible fluids. The purpose of this paper is
to understand whether compressible vortex sheets in three
dimensions, which are unstable in the regime of pure gas dynamics,
become stable under the magnetic effect in three-dimensional
magnetohydrodynamics (MHD). One of the main features is that the
stability problem is equivalent to a free boundary problem whose
free boundary is a characteristic surface,
%with the characteristic boundary,
which is more delicate than noncharacteristic free boundary
problems. Another feature is that the linearized problem for
current-vortex sheets in MHD does not meet the uniform
Kreiss-Lopatinskii condition.
These features cause additional analytical difficulties and
especially prevent a direct use of the standard Picard iteration to
the nonlinear problem. In this paper, we develop a nonlinear
approach to deal with these difficulties in
three-dimensional MHD. %In particular,
We first carefully formulate the linearized problem for the
current-vortex sheets to show rigorously that the magnetic effect
makes the problem weakly stable and establish energy estimates,
especially high-order energy estimates, in terms of the
nonhomogeneous terms and variable coefficients without loss of the
order. Then we exploit these results to develop a suitable iteration
scheme of Nash-Moser-H\"{o}rmander type and establish its
convergence,
%toward acompressible current-vortex sheet to
which leads to the existence and stability of compressible
current-vortex sheets, locally in time, in the three-dimensional
MHD.
\end{abstract}
\maketitle

\section{Introduction}

We are concerned with the existence and stability of compressible
current-vortex sheets in three-dimensional magnetohydrodynamics
(MHD). The motion of inviscid MHD fluids is governed by the
following system:
\begin{equation}\label{f1.1}
\begin{cases}
\partial_t\rho+ \nabla\cdot (\rho v)=0,\cr
\partial_t(\rho v)+ \nabla\cdot(\rho v\otimes v-H\otimes H)
  +\nabla(p+\frac{1}{2}|H|^2)=0,\cr
\partial_t H-\nabla\times (v\times H)=0,\cr
\partial_t\big(\rho(\frac{1}{2}(|{v}|^2+\frac{|H|^2}{\rho})+e )\big)
   + \nabla\cdot \big(\rho{v}(\frac{1}{2}|{v}|^2+e+\frac{p}{\rho})
    +H\times({v}\times H)\big)=0,
\end{cases}
\end{equation}
and
\begin{equation}\label{f1.2}
\nabla\cdot H=0,
\end{equation}
where $\rho, v=(v_1,v_2,v_3), H=(H_1,H_2,H_3),$ and $p$ are the
density, velocity, magnetic field, and pressure, respectively;
$e=e(\rho, S)$ is the internal energy; and $S$ is the entropy. The
Gibbs relation
$$
\theta dS=de -\frac{p}{\rho^2}d\rho
$$
implies the constitutive relations:
$$
p=\rho^2 e_\rho(\rho,S), \qquad \theta=e_S(\rho,S),
$$
where $\theta$ is the temperature. The quantity
$c=\sqrt{p_\rho(\rho, S)}$ is the sonic speed of the fluid.

For smooth solutions, the equations in (\ref{f1.1}) are equivalent
to
\begin{equation}\label{f1.3}
\begin{cases}
(\partial_t+ {v}\cdot\nabla)p+ \rho c^2 \nabla\cdot {v}=0,\cr \rho
(\partial_t+{v}\cdot\nabla) {v}+\nabla p- (\nabla\times H)\times
H=0,\cr (\partial_t+{v}\cdot\nabla)H-(H\cdot\nabla){v}+H
\nabla\cdot{v}=0,\cr (\partial_t+{v}\cdot\nabla)S=0,
\end{cases}
\end{equation}
which can be written as a $8\times 8$
symmetric hyperbolic system for
$U=(p,v,H,S)$ of the form: %\stepcounter{line}
\begin{equation}\label{f1.4}
B_0(U)\partial_t U+\sum_{j=1}^3B_j(U)\partial_{x_j}U=0.
\end{equation}

\bigskip
Compressible vortex sheets occur ubiquitously in nature (cf.
\cite{AM,BBR,CF,GM,SNW,VD,VN}) and are fundamental waves
in entropy solutions to the multidimensional hyperbolic systems of
conservation laws. For example, the vortex sheets are one of the
core waves in the Mach reflection configurations when a planar shock
hits a wedge and, more generally, in the two-dimensional Riemann
solutions, which are building blocks of general entropy solutions to
the Euler equations in gas dynamics. Therefore, understanding the
existence and stability of compressible vortex sheets is an
important step towards our full understanding of fluid motions and
the behavior of entropy solutions to the multidimensional hyperbolic
systems of conservation laws,
% influid dynamics and MHD,
along with the existence and stability of shock and rarefaction
waves (cf. Majda \cite{majda}, Glimm-Majda \cite{GM}, Alinhac
\cite{ali}, and M\'{e}tivier \cite{met2}; also see
\cite{met1,GMWZ}).

For the Euler equations in two-dimensional gas dynamics, the
classical linearized stability analysis on compressible vortex
sheets predicts stability when the Mach number $M>\sqrt{2}$ since
there are no growing modes in this case, and instability when
$M<\sqrt{2}$ since there are exponentially exploding modes of
instability (see Miles \cite{miles} for the definitive treatment).
The local nonlinear stability for the two-dimensional case with
$M>\sqrt{2}$ was recently established by Coulombel-Secchi
\cite{cou1,cou2} when the initial data is a small perturbation of a
planar vortex sheet. In a series of papers, Artola and Majda
\cite{AM} studied the stability of vortex sheets by using the
argument of nonlinear geometric optics to analyze the interaction
between vortex sheets and highly oscillatory waves. They first
observed the generation of three distinct families of kink modes
traveling along the slip stream bracketed by shocks and rarefaction
waves and provided a detailed explanation of the instability of
supersonic vortex sheets even when $M>\sqrt{2}$. Therefore, one can
not expect
the nonlinear stability globally for
two-dimensional compressible vortex sheets.

For the Euler equations in three-dimensional gas dynamics, the
situation is even more complicated. In fact, it is well known that
every compressible vortex sheet is violently unstable and this
violent instability is the analogue of the Kelvin-Helmholtz
instability for incompressible fluids.

The purpose of this paper is to understand whether compressible
vortex sheets in three dimensions, which are unstable in the regime
of pure gas dynamics, become stable under the magnetic effect, that
is, the nonlinear stability of current-vortex sheets in
three-dimensional MHD. One of the main features is that the
stability problem is equivalent to a free boundary problem whose
free boundary is a characteristic surface,
%with the characteristic boundary,
which is more delicate than the noncharacteristic free boundary
problems that are endowed with a strict jump of the normal velocity
on the free boundary (cf. \cite{majda,met1,met2}; also
\cite{CFe,CF}). Another feature is that the linearized problem for
current-vortex sheets in MHD does not meet the uniform
Kreiss-Lopatinskii condition \cite{kreiss,majda}; see the analysis
in Blokhin-Trakhinin \cite{blo} and Trakhinin \cite{tra}.
These features cause additional analytical difficulties and,
in particular, they prevent a direct use of the standard Picard
iteration to prove the existence of solutions to the nonlinear
problem.

In this paper, we develop a nonlinear approach to deal with these
difficulties in three-dimensional MHD, based on the previous results
mentioned above. We first carefully formulate the linearized problem
for the current-vortex sheet to show rigorously that, although any
compressible vortex sheet may be violently unstable in the regime of
pure gas dynamics and does not meet the uniform Kreiss-Lopatinskii
condition, the magnetic effect makes the linearized problem
weakly stable as observed in \cite{blo,tra}. In
particular, we successfully establish high-order energy estimates of
the solutions for the linearized problem in terms of the
nonhomogeneous terms and variable coefficients without loss of the
order. Then we exploit these results to develop a suitable iteration
scheme of Nash-Moser-H\"{o}rmander type and establish its
convergence,
which leads to the existence and stability of compressible
current-vortex sheets, locally in time, in the three-dimensional
MHD.

We remark that, in order to establish the energy estimates,
especially {\it high-order energy estimates}, of solutions  to the
linearized problem derived from the current-vortex sheet, we have
identified a well-structured decoupled formulation so that the
linear problem \eqref{f3.1}--\eqref{f3.3} is decoupled into one
standard initial-boundary value problem \eqref{f3.10} for a
symmetric hyperbolic equations and another problem \eqref{f3.11} for
an ordinary differential equation for the front. This decoupled
formulation is much more convenient and simpler than that in
(58)--(60) given in \cite{tra} and is essential for us to establish
the desired high-order energy estimates of solutions, which is one
of the key ingredients for developing the suitable iteration scheme
of Nash-Moser-H\"{o}rmander type that converges.

This paper is organized as follows. In Section 2, we first set up
the current-vortex sheet problem as a free boundary problem, and
then we reformulate this problem into a fixed initial-boundary value
problem and state the main theorem of this paper. In Sections 3--4,
we first carefully formulate the linearized problem for compressible
current-vortex sheets and identify the decoupled formulation; and
then we make energy estimates, especially high-order energy
estimates, of the solutions in terms of the nonhomogeneous terms and
variable coefficients without loss of the order. Then, in Section 5,
we analyze the compatibility conditions and construct the zero-th
order approximate solutions for the nonlinear problem, which is a
basis for our iteration scheme. In Section 6, we develop a suitable
iteration scheme of Nash-Moser-H\"{o}rmander type. In Section 7, we
establish the convergence of the iteration scheme towards a
compressible current-vortex sheet. Finally, in Section 8, we
complete several necessary estimates of the iteration scheme used in
Section 7.

\section{Current-Vortex Sheets and Main Theorem}

In this section, we first set up the current-vortex sheet problem as
a free boundary problem and then reformulate this problem into a
fixed initial-boundary value problem, and finally we state the main
theorem of this paper.

Let a piecewise smooth function $U(t,x)$ be an entropy solution to
(\ref{f1.1}) with the form:
\begin{equation}\label{f2.1}
U(t,x)=\begin{cases}
U^+(t,x) \quad \mbox{when}\,\, x_1>\psi(t,x_2,x_3),\\
U^-(t,x) \quad \mbox{when}\,\, x_1<\psi(t,x_2,x_3). \end{cases}
\end{equation}
Then, on the discontinuity front $\Gamma=\{x_1=\psi(t,x_2,x_3)\}$,
$U(t,x)$ must satisfy the Rankine-Hugoniot conditions:
%\stepcounter{line}
\begin{equation}\label{f2.2}
\begin{cases}
[m_N]=0,\qquad [H_N]=0,\cr m_N[v_N]+[q]=0, \qquad
m_N[v_\tau]=H_N[H_\tau],\cr m_N[\frac{H_k}{\rho}]=H_N[v_k]\qquad
{\rm for}~k=1,2,3,\cr
m_N[e+\frac{1}{2}(|v|^2+\frac{|H|^2}{\rho})]+[qv_N-H_N(H\cdot v)]=0,
\end{cases}
\end{equation}
where $[\alpha]$ denotes the jump of the function $\alpha$ on the
front $\Gamma$, $(v_N, v_\tau)$ (resp. $(H_N, H_\tau)$) are the
normal and tangential components of $v$  (resp. $H$) on $\Gamma$,
i.e.,
\begin{eqnarray*}
&&v_N:=v_1-\psi_{x_2}v_2-\psi_{x_3}v_3, \qquad
v_\tau=(v_{\tau_1},v_{\tau_2})^\top:=(\psi_{x_2}v_1+v_2,\psi_{x_3}v_1+v_3)^\top,\\
&&H_N:=H_1-\psi_{x_2}H_2-\psi_{x_3}H_3, \qquad
H_\tau=(H_{\tau_1},H_{\tau_2})^\top:=(\psi_{x_2}H_1+H_2,\psi_{x_3}H_1+H_3)^\top,
\end{eqnarray*}
$m_N=\rho(v_N-\psi_t)$ is the mass transfer flux, and
$q=p+\frac{|H|^2}{2}$ is the total pressure.

\medskip
If $m_N=0$ on $\Gamma$, then $(U^\pm,\Gamma)$ is a contact
discontinuity for (\ref{f1.1}). In this paper, we focus on the case:
\begin{eqnarray}
&&H^+_N=H^-_N=0, \label{f2.3}\\
&&H^+_\tau\not\parallel H^-_\tau, \label{f2.3a}
\end{eqnarray}
for which $U(t,x)$ is called a compressible current-vortex sheet.
Then the current-vortex sheet is determined by \eqref{f2.3} and
\begin{equation}\label{f2.4}
\psi_t=v_N^\pm, \qquad [p+\frac{|H|^2}{2}]=0,
\end{equation}
on which generically
\begin{equation}\label{f2.5}
([\rho], [v_\tau], [S])\neq 0.
\end{equation}
When $H\equiv 0$, $U(t,x)$ becomes a classical compressible vortex
sheet in fluid mechanics without magnetic effect.

\medskip
Since
$H_N^+=H^-_N$ on $\Gamma=\{x_1=\psi(t,x_2,x_3)\}$, it is easy to see
that {\it condition \eqref{f2.3a} is equivalent to the condition:
\begin{equation}\label{f2.6}
\left(
\begin{array}{l}
H_2^+\\
H_3^+\end{array}\right)\not\parallel
\left(
\begin{array}{l}
H_2^-\\
H_3^-\end{array}\right)\qquad{\rm on}\,\,\Gamma.
\end{equation}}
Indeed, when $H_N^+=H^-_N$ on $\Gamma$, we have
$$
H_\tau^\pm=((1+\psi_{x_2}^2)H_2^\pm+\psi_{x_2}\psi_{x_3}H_3^\pm,
\psi_{x_2}\psi_{x_3}H_2^\pm+(1+\psi_{x_2}^2)H_3^\pm)^\top
$$
which implies that $\alpha H_\tau^++\beta H_\tau^-=0$ if and only if
$\alpha H_k^++\beta H_k^-=0$ hold for $k=2,3$.

\medskip
Furthermore, system (\ref{f1.1}) has the following property: {\it
Let $(U^\pm, \psi)$ be the current-vortex sheet to \eqref{f1.1} for
$t\in [0,T)$ for some $T>0$. Then, if $ \nabla\cdot H^\pm(0,x)=0, $
we have
\begin{equation}\label{f2.7a}
 \nabla\cdot H^\pm(t,x)=0\qquad
\mbox{for all}\,\, t\in [0, T).
\end{equation}}
This can be seen as follows: From (\ref{f1.3}), we find that, on
both sides of $\Gamma$,
\begin{equation}\label{f2.7}
(\partial_t + v^\pm\cdot\nabla)(\nabla\cdot H^\pm)+ (\nabla\cdot
{v}^\pm)(\nabla\cdot H^\pm)=0.
\end{equation}
On $\Gamma$, we use $\psi_t= v^\pm_N$ to obtain
$$
\partial_t +{v}^\pm\cdot\nabla
= \partial_t+(v^\pm_N+\psi_{x_2}v^\pm_2+\psi_{x_3}v^\pm_3)
 \partial_{x_1}+v^\pm_2\partial_{x_2}
  +v^\pm_3\partial_{x_3}=\tau^\pm\cdot(\partial_t, \nabla),
$$
where the vectors $\tau^\pm=(1,
v^\pm_N+\psi_{x_2}v^\pm_2+\psi_{x_3}v^\pm_3, v^\pm_2, v^\pm_3)^\top$
are orthogonal to the space-time normal vector ${n}=(\psi_t, -1,
\psi_{x_2}, \psi_{x_3})^\top$ on $\Gamma$. Thus, $\partial_t+
v^\pm\cdot\nabla$ are tangential derivative operators to $\Gamma$.
Thus, equations in (\ref{f2.7}) are transport equations for
$\nabla\cdot H$ on both sides of $\Gamma$, which implies
\eqref{f2.7a}.

\bigskip Set
$$
D(\lambda, U):=\left(
\begin{array}{cc}
\tilde{D}(\lambda, U) & 0\\
0 & 1\end{array}\right)
\qquad
\mbox{with}
\quad \tilde{D}(\lambda,
U):=\left(
\begin{array}{ccc}
1 & \frac{\lambda}{\rho c^2}H^\top & O_{1\times 3}\\
\lambda\rho H & I_3 & -\rho\lambda I_3\\
0_{3\times 1} & -\lambda I_3 & I_3\end{array}\right).
$$
As in \cite{tra}, property \eqref{f2.7a} implies that system
(\ref{f1.1})--(\ref{f1.2}) is equivalent to the following system on
both sides of $\Gamma$:
%\stepcounter{line}
\begin{equation}\label{f2.8}
D(\lambda^\pm, U^\pm)\Big(B_0(U^\pm)\partial_tU^\pm +\sum_{j=1}^3
B_j(U^\pm) \partial_{x_j}U^\pm\Big) +\lambda^\pm {G}^\pm \,
\nabla\cdot H^\pm=0,
\end{equation}
provided $\nabla\cdot H^\pm(0,x)=0$, where
$$
{G}^\pm=-(1, 0,0,0,H^\pm, 0)^\top,
$$
and $\lambda^\pm=\lambda^\pm(v_2^\pm, v_3^\pm, H_2^\pm, H_3^\pm)$
will be determined later so that
\begin{equation}\label{f2.9}
(\lambda^\pm)^2<\frac{(c^\pm)^2}{\rho^\pm
(c^\pm)^2+|H^\pm|^2}
\end{equation}
are satisfied for the sound speeds $c^\pm$.

We rewrite system (\ref{f2.8}) as
\begin{equation}\label{f2.10}
A_0(U^\pm)\partial_t
U^\pm+\sum_{j=1}^3A_j(U^\pm)\partial_{x_j}U^\pm=0,
\end{equation}
which is also symmetric hyperbolic. Then the problem of existence
and stability of current-vortex sheets of the form (\ref{f2.1}) can
be formulated as the following free boundary problem:

\medskip
{\bf Free Boundary Problem}:\, {\it Determine $U^\pm(t,x)$ and a
free boundary $\Gamma=\{x_1=\psi(t,x_2,x_3)\}$ for $t>0$ such that
\begin{equation}\label{f2.11}
\begin{cases}
A_0(U^\pm)\partial_t
U^\pm+\sum_{j=1}^3A_j(U^\pm)\partial_{x_j}U^\pm=0
\qquad\mbox{for}\,\, \pm(x_1-\psi(t,x_2, x_3))>0,\\
U|_{t=0}=\begin{cases}U_0^+(x)
\quad \mbox{when}\,\, x_1>\psi_0(x_2, x_3),\\
U_0^-(x) \quad\mbox{when}\,\, x_1<\psi_0(x_2, x_3)
\end{cases}
\end{cases}
\end{equation}
satisfying the jump conditions on $\Gamma$:
\begin{equation}\label{f2.12}
\psi_t=v_N^\pm, \quad H_N^\pm=0, \quad [p+\frac{|H|^2}{2}]=0, \quad
H^+_\tau\not\parallel H^-_\tau,
\end{equation}
where $\psi_0=\psi(0,x_1,x_2)$, and condition {\rm (\ref{f2.9})} is
satisfied for $\lambda^\pm$ given by {\rm (\ref{f3.9})} later.}

\medskip
For a piecewise smooth solution $U(t,x)$ of problem
(\ref{f2.11})--(\ref{f2.12}) with smooth $(U^\pm, \psi)$ for $t\in
[0,T)$, it can be similarly verified that, if $\nabla\cdot
H^\pm(0,x)=0$, then $\nabla\cdot H^\pm(t,x)=0$ for all $t\in [0,
T).$  Hence the equations in (\ref{f2.11}) and
(\ref{f1.1})--(\ref{f1.2}) are equivalent for such a solution.

Furthermore, unlike the shock case, the new difficulty here is that
the free boundary $\Gamma$ in problem (\ref{f2.11})--(\ref{f2.12})
is now a characteristic surface of system (\ref{f2.11}),
rather than the noncharacteristic boundary surface endowed with a
strict jump of the normal velocity as in the shock case (see
\cite{CFe,CF,met1,majda,met2}).

\medskip
To deal with such a free boundary problem,
it is convenient to employ the standard partial hodograph
transformation:
\begin{equation}\label{f2.13}
\begin{cases}t={\tilde t}, \quad x_2={\tilde x_2}, \quad x_3={\tilde
x_3},\cr x_1=\Psi^\pm({\tilde t}, {\tilde x_1}, {\tilde x_2},
{\tilde x_3})
\end{cases}
\end{equation}
with $\Psi^\pm$ satisfying
%\stepcounter{line}
\begin{equation}\label{f2.14}
\begin{cases}
\pm\partial_{\tilde x_1}\Psi^\pm\ge  \kappa>0,\cr \Psi^+|_{{\tilde
x_1}=0}=\Psi^-|_{{\tilde x_1}=0}=\psi({\tilde t}, {\tilde x_2},
{\tilde x_3})
\end{cases}
\end{equation}
for some constant $\kappa>0$. Under (\ref{f2.13}), the domains
$\Omega^\pm=\{\pm(x_1-\psi(t,x_2,x_3))>0\}$ are transformed into
$\{{\tilde x_1}>0\}$ and the free boundary $\Gamma$ into the fixed
boundary $\{\tilde{x}_1=0\}$.

Then we define ${\tilde U}^\pm({\tilde t}, {\tilde x}):=U^\pm(t,x)$.
{}From the first jump condition in (\ref{f2.12}), the natural
candidates of $\Psi^\pm$ are those that satisfy the following
eikonal equations:
\begin{equation}\label{f2.15}
\partial_t \Psi^\pm-v_1^\pm+v_2^\pm\partial_{x_2}\Psi^\pm
+v_3^\pm\partial_{x_3}\Psi^\pm=0 \qquad\mbox{in}\,\,\{x_1>0\},
\end{equation}
where we have dropped the tildes in the formula for simplicity of
notations.

It is easy to check that, under (\ref{f2.13}), ${\tilde
U}^\pm({\tilde t}, {\tilde x})=U^\pm(t,x)$ satisfy
%\stepcounter{line}
\begin{equation}\label{f2.16}
L(U^\pm, \Psi^\pm)U^\pm=0 \qquad {\rm in}~\{x_1>0\},
\end{equation}
the boundary conditions on $\{x_1=0\}$:
\begin{equation}\label{f2.17}
\begin{cases}
\Psi^+|_{x_1=0}=\Psi^-|_{x_1=0}=\psi,\cr B(U^+, U^-,
\psi)|_{x_1=0}=0,
\end{cases}
\end{equation}
and the initial condition at $\{t=0\}$:
\begin{equation}\label{f2.18}
(U^\pm, \psi)|_{t=0}=(U_0^\pm(x), \psi_0(x_2, x_3)),
\end{equation}
where the tildes have also been dropped for simplicity of notations,
$$
L(U, \Psi)V=A_0(U)\partial_t V +{\bar A}_1(U, \Psi)
\partial_{x_1}V+\sum_{j=2}^3A_j(U)
\partial_{x_j}V
$$
with ${\bar A}_1(U, \Psi)=
\frac{1}{\partial_{x_1}\Psi}\big(A_1(U)-\partial_t\Psi
A_0(U)-\sum_{j=2}^3\partial_{x_j}\Psi A_j(U)\big)$,
$$
B(U^+, U^-, \psi)=\left(\partial_t\psi- U_{v,N}^\pm,
U_{H,N}^\pm, q^+-q^-\right)^\top
$$
with $U_{v,N}^\pm=U_2^\pm-\psi_{x_2}U_3^\pm-\psi_{x_3}U_4^\pm,
U_{H,N}^\pm=U_5^\pm-\psi_{x_2}U_6^\pm-\psi_{x_3}U_7^\pm,$ and
$$
q=U_1+\frac{1}{2}|U_H|^2 \qquad\mbox{for}\,\, U_H=(U_5, U_6,
U_7)^\top.
$$
With these, the free boundary problem has been reduced into the
fixed {\bf initial-boundary value problem
\eqref{f2.15}--\eqref{f2.18}}.

\medskip
To solve \eqref{f2.15}--\eqref{f2.18}, as in \cite{ali, che, gue},
it is natural to introduce the tangential vector $M=(M_1,M_2, M_3)$
of $\{x_1=0\}$:
$$
M_1=\sigma(x_1)\partial_{x_1} \quad M_2=\partial_{x_2}, \quad
M_3=\partial_{x_3},
$$
with
$$
\sigma(x_1):=
\begin{cases}
x_1 \quad \mbox{when}\,\, 0\le x_1\le 1,\\
2 \quad \,\,\,\mbox{when}\,\, x_1\ge 2,
\end{cases}
$$
and the weighted Sobolev spaces defined on $\Omega_T:=\{(t,x)\in [0,
T]\times \R^3\, :\, ~x_1>0\}$:
$$
B^s_\mu(\Omega_T):=\{u\in L^2(\Omega_T)\,:\, e^{-\mu
t}M^\alpha\partial_{x_1}^k u\in  L^2(\Omega_T)\quad\mbox{for}\,
|\alpha|+2k\le s\}
$$
for all $s\in \N$ and $\mu>0$, imposed with the
norms
$$
\|u\|_{s,\mu,T}:=\big(\int_0^T\|u(t,\cdot)\|_{s,\mu}^2
dt\big)^{1/2},
$$
where
$$ \|u(t,\cdot)\|_{s,\mu}^2:=\sum_{|\alpha|+2k\le
s}\mu^{2(s-|\alpha|-2k)} \| e^{-\mu t}M^\alpha \partial_{x_1}^k u(t,
\cdot)\|_{L^2}^2.
$$
We will also use similar notation as above for the spaces with
$\mu=0$, $B^s(\Omega_T)$, whose norm is defined by
$$
\|u\|_{s,T}:=\big(\sum_{|\alpha|+2k\le s} \| M^\alpha
\partial_{x_1}^k u\|_{L^2(\Omega_T)}^2\big)^{1/2}.
$$
Also denote $b\Omega_T:=\{(t,x_2,x_3)\, :\, t\in [0, T],\, (x_2,
x_3)\in \R^2\}$.

Then, when the initial data $(U_0^\pm, \psi_0)$ is a small
perturbation of a plane current-vortex sheet
$(\bar{U}^\pm,\bar{\psi})$ for constant states $\bar{U}^\pm$ and
$\bar{\psi}=0$ satisfying \eqref{f2.3}--\eqref{f2.4}, we have the
following main theorem of this paper.

\bigskip
{\bf Theorem 2.1 (Main Theorem).} {\it Assume that, for any fixed
$\alpha\ge 14$ and $s\in[\alpha+5, 2\alpha-9]$,  the initial data
functions $\psi_0\in H^{2s+3}(\R^2)$ and $U_0^\pm-\bar{U}^\pm\in
B^{2(s+2)}(\R^3_+)$ satisfy the compatibility conditions of problem
\eqref{f2.15}--\eqref{f2.18} up to order $s+2$.
Then there exists a solution $(U^\pm, \Psi^\pm)$ of the
initial-boundary value problem \eqref{f2.15}--\eqref{f2.18} such
that
$$ (U^\pm-U_a^\pm, \Psi^\pm-\Psi_a^\pm)\in B^\alpha(\Omega_T)
\quad\mbox{and}\quad \psi=\Psi^\pm|_{x_1=0}\in
H^{\alpha-1}(b\Omega_T)
$$
for some functions $U_a^\pm$ and $\Psi_a^\pm$ satisfying
$$
U_a^\pm-\bar{U}^\pm, \Psi_a^\pm\mp x_1 \in
B^{s+3}(\R_+\times\R^3_+).
$$
}

We remark that, since $\pm\partial_{\tilde x_1}\Psi^\pm\ge
\kappa>0$, the corresponding vector function of the solution
$(U^\pm, \Psi^\pm)$ in Theorem 2.1 under the inverse of the partial
hodograph transform is a solution of the free boundary problem, i.e.
the current-vortex sheet problem, which implies the existence and
stability of compressible current-vortex sheets to
\eqref{f2.1}--\eqref{f2.2} under the initial perturbation. To
achieve this, we start with the linear stability in Sections 3--4
and then establish the nonlinear stability by developing nonlinear
techniques in Sections 5--7. Some estimates used in Sections 5--7
for the iteration scheme of Nash-Moser-H\"{o}rmander type are given
in Section 8.

\section{Linear Stability I: Linearized Problem}

To study the linear stability of current-vortex sheets, we first
derive a linearized problem from the nonlinear problem
\eqref{f2.15}--\eqref{f2.18}.
%As in Proposition 1.3.1 of \cite{met2},
By a direct calculation, we have
$$
\frac{d}{ds}\big(L(U+sV, \Psi+s\Phi)(U+sV)\big)|_{s=0} =L(U,
\Psi)W+E(U, \Psi)W+\frac{\Phi}{\partial_{x_1}\Psi}\partial_{x_1}
(L(U, \Psi)U),
$$
where
$$
W=V-\frac{\Phi}{\partial_{x_1}\Psi}\partial_{x_1}U
$$
is the good unknown as introduced in \cite{ali} (see also
\cite{met1, met2}) and
$$
E(U, \Psi)W =W\cdot \nabla_U({\bar
A}_1(U,\Psi))\partial_{x_1}U+\sum_{j=2}^3W\cdot \nabla_UA_j(U)
\partial_{x_j}U.
$$

Then we obtain the following linearized problem of
\eqref{f2.15}--(\ref{f2.18}):
\begin{equation}\label{f3.1}
L(U^\pm, \Psi^\pm)W^\pm+E(U^\pm, \Psi^\pm)W^\pm=F^\pm\qquad
\mbox{in}\,\, \{x_1>0\}
\end{equation}
with the boundary conditions on $\{x_1=0\}$:
\begin{equation}\label{f3.2}
\begin{cases}
\phi_t-(W_2^\pm-\psi_{x_2}W_3^\pm-\psi_{x_3}W_4^\pm)+U_3^\pm
\phi_{x_2}+U_4^\pm \phi_{x_3}=h_1^\pm,\cr
W_5^\pm-\psi_{x_2}W_6^\pm-\psi_{x_3}W_7^\pm -U_6^\pm
\phi_{x_2}-U_7^\pm \phi_{x_3}=h_2^\pm,\cr
W^+_1-W^-_1+\sum_{j=5}^7(U_j^+W_j^+-U_j^-W_j^-)=h_3,
\end{cases}
\end{equation}
and the initial condition at $\{t=0\}$:
\begin{equation}\label{f3.3}
%\label{l-init}
(W^\pm, \phi)|_{t=0}=0,
\end{equation}
for some functions $F^\pm$ and $h:=(h_1^\pm, h_2^\pm, h_3)^\top$.

\medskip
To separate the characteristic and noncharacteristic components of
the unknown $W^\pm$, we introduce $J^\pm=J(U^\pm, \Psi^\pm)$ as a
$8\times 8$ regular matrix such that
\begin{equation}\label{f3.4}
X^\pm=(J^\pm)^{-1}W^\pm
\end{equation}
satisfies
$$
\begin{cases}
X_1^\pm=W_1^\pm+\sum_{j=5}^7U_j^\pm W_j^\pm,\cr
X_2^\pm=W_2^\pm-\Psi^\pm_{x_2}W_3^\pm-\Psi^\pm_{x_3}W_4^\pm,\cr
X_5^\pm=W_5^\pm-\Psi^\pm_{x_2}W_6^\pm-\Psi^\pm_{x_3}W_7^\pm,\cr
(X_3^\pm, X_4^\pm, X_6^\pm, X_7^\pm, X_8^\pm)=(W_3^\pm, W_4^\pm,
W_6^\pm, W_7^\pm, W_8^\pm).
\end{cases}
$$
Then, under transformation (\ref{f3.4}), problem
(\ref{f3.1})--\eqref{f3.3} for $(W^\pm, \phi)$ is equivalent to
the following problem for $(X^\pm,\phi)$:
%\stepcounter{line}
\begin{equation}\label{f3.5}
{\tilde L}(U^\pm, \Psi^\pm)X^\pm+{\tilde E}(U^\pm,
\Psi^\pm)X^\pm={\tilde F}^\pm\qquad \mbox{in}\,\, \{x_1>0\}
\end{equation}
with the boundary conditions on $\{x_1=0\}$:
\begin{equation}\label{f3.6}
\begin{cases}
\phi_t-X_2^\pm+U_3^\pm \phi_{x_2}+U_4^\pm \phi_{x_3}=h_1^\pm,\cr
X_5^\pm-U_6^\pm \phi_{x_2}-U_7^\pm \phi_{x_3}=h_2^\pm,\cr
X^+_1-X^-_1=h_3,
\end{cases}
\end{equation}
and the initial condition at $\{t=0\}$:
\begin{equation}\label{f3.7}
(X^\pm, \phi)|_{t=0}=0,
%\phi_0(x_2,x_3),
\end{equation}
where ${\tilde F}^\pm=(J^\pm)^\top F^\pm$,
$$
{\tilde L}(U^\pm, \Psi^\pm)
={\tilde A}_0(U^\pm, \Psi^\pm)\partial_t
+\sum_{j=1}^3{\tilde A}_j(U^\pm, \Psi^\pm)\partial_{x_j}
$$
with ${\tilde A}_1(U^\pm, \Psi^\pm)=(J^\pm)^\top {\bar A}_1(U^\pm,
\Psi^\pm) J^\pm$, ${\tilde A}_j(U^\pm, \Psi^\pm)=(J^\pm)^\top
A_j(U^\pm) J^\pm$, $j\neq 1$, and
$$
{\tilde E}(U^\pm, \Psi^\pm)X^\pm =(J^\pm)^\top E(U^\pm,
\Psi^\pm)J^\pm X^\pm + (J^\pm)^\top (L(U^\pm, \Psi^\pm) J^\pm)
X^\pm.
$$

For simplicity of notations, we will drop the tildes in
(\ref{f3.5}). By a direct calculation, we see that $A_1(U^\pm,
\Psi^\pm)$ can be decomposed into three parts:
$$
A_1(U^\pm, \Psi^\pm)=A_1^{\pm, 0}+A_1^{\pm, 1}+A_1^{\pm, 2}
$$
with
$$
A_1^{\pm, 0}=\frac{1}{\partial_{x_1}\Psi^\pm}\left[
\begin{array}{cc}
0 & {a}\\
{a}^\top & O_{7\times 7}\end{array}\right], \quad A_1^{\pm,
1}=\frac{\partial_{t}\Psi^\pm-U^\pm_{v,N}}{\partial_{x_1}\Psi^\pm}
{\tilde A}_1^{\pm, 1}, \quad A_1^{\pm,
2}=\frac{U^\pm_{H,N}}{\partial_{x_1}\Psi^\pm} {\tilde A}_1^{\pm,
2}
$$
where ${a}=(1, 0,0,-\lambda^\pm, 0, 0, 0)$, $U_{v,N}^\pm=U_2^\pm-
\Psi_{x_2}U_3^\pm- \Psi_{x_3}U_4^\pm$, and $U_{H, N}^\pm=U_5^\pm-
\Psi_{x_2}U_6^\pm- \Psi_{x_3}U_7^\pm$.

When the state $(U^\pm, \Psi^\pm)$ satisfies the boundary conditions
given in (\ref{f2.17}), we know from (\ref{f3.5}) that, on
$\{x_1=0\}$,
%\stepcounter{line}
\begin{equation}\label{f3.8}
<\left(\begin{array}{cc}
A_1(U^+, \Psi^+) & 0 \\ 0 &  A_1(U^-, \Psi^-)\end{array}\right)
\left(\begin{array}{c}
X^+ \\X^-\end{array}\right), \left(\begin{array}{c}
X^+ \\X^-\end{array}\right)>
=2X_1^\pm[X_2-\lambda X_5],
\end{equation}
when $[X_1]=0$ on $\{x_1=0\}$.

Furthermore, from (\ref{f3.6}), we have
$$
[X_2-\lambda X_5]
=\phi_{x_2}[U_3-\lambda U_6]+\phi_{x_3}[U_4-\lambda U_7]-[h_1+\lambda h_2].
$$
Thus, from assumption \eqref{f2.3a} that is equivalent to
\eqref{f2.6},
there exists a unique $\lambda^\pm=\lambda^\pm(v_2^\pm, v_3^\pm,
H_2^\pm, H^\pm_3)$ such that
\begin{equation}\label{f3.9}
\left(\begin{array}{c} v_2^+
\\v_3^+\end{array}\right)-\left(\begin{array}{c} v_2^-
\\v_3^-\end{array}\right)=\lambda^+ \left(\begin{array}{c} H_2^+
\\H_3^+\end{array}\right) -\lambda^- \left(\begin{array}{c} H_2^-
\\H_3^-\end{array}\right),
\end{equation}
that is, $[U_3-\lambda U_6]=[U_4-\lambda U_7]=0$, which implies
$$
[X_2-\lambda X_5]=-[h_1+\lambda h_2].
$$

Therefore, with the choice of $\lambda^\pm$ in \eqref{f3.9}, problem
(\ref{f3.5})--\eqref{f3.7} is decomposed into
\begin{equation}\label{f3.10}
\begin{cases} L(U^\pm, \Psi^\pm)X^\pm+E(U^\pm,
\Psi^\pm)X^\pm=F^\pm\qquad \mbox{in}\,\, \{x_1>0\},\cr [X_2-\lambda
X_5]=-[h_1+\lambda h_2] \qquad\mbox{on}\,\,\{x_1=0\},\cr
X^+_1-X^-_1=h_3 \qquad\mbox{on}\,\, \{x_1=0\},\cr X^\pm|_{t\le 0}=0,
\end{cases}
\end{equation}
which is maximally dissipative in the sense of Lax-Friedrichs
\cite{LF,majda}, and
\begin{equation}\label{f3.11}
\begin{cases}\phi_t=X_2^\pm-(U_3^\pm, U^\pm_4)
\left(\begin{array}{cc} U_6^+ & U_7^+ \\ U_6^- &
U_7^-\end{array}\right)^{-1} \left(\begin{array}{c} X_5^+-h_2^+ \\
X_5^--h_2^-
\end{array}\right)+h_1^\pm\qquad\mbox{on}\,\, \{x_1=0\},\cr
\phi|_{t=0}=0.
\end{cases}
\end{equation}

\medskip
{\bf Remark 3.1.} To have the identities in (\ref{f3.11}), it
requires the compatibility, that is, the right-hand side of
(\ref{f3.11}) for the plus sign is equal to the term for the minus
sign. This is also guaranteed by the choice of $\lambda^\pm$ in
(\ref{f3.9}).

\section{Linear Stability II: Energy Estimates}

In this section, we establish the energy estimates for the linear
problem (\ref{f3.1})--(\ref{f3.3}), which implies the existence and
uniqueness of solutions. We know from Section 3 that it suffices to
study the {\bf linear problem (\ref{f3.10})--(\ref{f3.11})}.

As noted as above, when the state $(U^\pm, \Psi^\pm)$ satisfies the
eikonal equation (\ref{f2.15}) in $\{x_1>0\}$ and the boundary
conditions on $\{x_1=0\}$ in (\ref{f2.17}), we know that the
boundary $\{x_1=0\}$ is a uniform characteristic of multiplicities
$12$ of the equations in (\ref{f3.10}).

First, we have the following elementary properties in the space
$B^s(\Omega_T)$, whose proof can be found in \cite{ali}.

\vspace{.1in} {\bf Lemma 4.1.} {\it {\rm (i)} For any fixed $s>
5/2$, the identity mapping is a bounded embedding from
$B^s(\Omega_T)$ to $L^\infty(\Omega_T)$;

{\rm (ii)} If $s> 1$ and $u\in B^s(\Omega_T)$, then $u|_{x_1=0}\in
H^{s-1}(b\Omega_T)$, the classical Sobolev space. Conversely, if
$v\in H^\sigma(b\Omega_T)$ for a fixed $\sigma>0$, then there exists
$u\in B^{\sigma+1}(\Omega_T)$ such that $u|_{x_1=0}=v$.}

\vspace{.1in} Denote by $coef (t,x)$ all the coefficient functions
appeared in the equation in (\ref{f3.10}), and
$\dot{coef}(t,x)={coef}(t,x)-{coef}(0)$. With the decoupled
formulation (\ref{f3.10})--(\ref{f3.11}) for the linearized problem,
we can now employ the Lax-Friedrichs theory \cite{LF,majda} to
establish the desired energy estimates, especially the high-order
energy estimates, of solutions for the linear problem.

\vspace{.1in} {\bf Theorem 4.1.} {\it For any fixed $s_0>17/2$,
there exist constants $C_0$ and $\mu_0$ depending only on
$\|\dot{coef}\|_{s_0,T}$ for the coefficient functions in
\eqref{f3.10} such that, for any $s\ge s_0$ and $\mu\ge \mu_0$, the
estimate
\begin{equation}\label{f4.1}
\begin{aligned}
&\max_{0\le t\le T}\|X^\pm (t)\|_{s, \mu}^2 +\mu \|X^\pm\|_{s, \mu,
T}^2 \\
&\le \frac{C_0}{\mu}\big(\|{F}^\pm\|_{s, \mu,
T}^2+\|{h}\|_{H^{s+1}_\mu(b\Omega_T)}^2+\|\dot{coef}\|_{s,
\mu, T}^2 (\|{F}^\pm\|_{s_0, T}^2
+\|{h}\|_{H^{s_0+1}(b\Omega_T)}^2)\big)
\end{aligned}
\end{equation}
holds provided that the eikonal equations {\rm (\ref{f2.15})} and
$$U_{H,N}^\pm:=U_5^\pm-\partial_{x_2}\Psi^\pm U_6^\pm-\partial_{x_3}\Psi^\pm U_7^\pm=0$$
are valid for $(U^\pm, \Psi^\pm)$ in a neighborhood of $\{x_1=0\}$,
where ${h}=(h_1^\pm, h_2^\pm, h_3)^\top$ and the norms  in
$H^s_\mu(b\Omega_T)$ are defined as that of $B_\mu^s(\Omega_T)$ with
functions independent of $x_1$.}

\vspace{.1in} {\it Proof:} The proof is divided into two steps.

{\it Step 1}. We first study problem (\ref{f3.10}) with homogeneous
boundary conditions:
%\stepcounter{line}
\begin{equation}\label{f4.2}
h_1^\pm=h_2^\pm=h_3=0 \qquad {\rm on}~\{x_1=0\}.
\end{equation}
Define
$$
P^\pm=\left(\begin{array}{ccc}
I_2 &  0 & 0\\
0 & P_1 & 0\\
0 & 0 & I_3\end{array}\right)\left(\begin{array}{cc}
P_2^\pm &  0\\
0 & I_5\end{array}\right)
$$
with
$$
P_1=\left(\begin{array}{ccc}
0 &  1 & 0\\
0 & 0 & 1\\
1 & 0 & 0\end{array}\right), \quad
P_2^\pm=\left(\begin{array}{ccc}
1 &  0 & 0\\
0 & 1 & \lambda^\pm\\
0 & 0 & 1\end{array}\right).
$$
Then
$$
{\mathcal A}_1^\pm=(P^\pm)^\top A_1^\pm P^\pm=
\left(\begin{array}{cc}
A_{11}^{\pm, 1} &  A_{12}^{\pm, 1}\\
A_{21}^{\pm, 1} &A_{22}^{\pm, 1}\end{array}\right)$$
with
%\stepcounter{line}
\begin{equation}\label{f4.3}
A_{11}^{\pm, 1}=
\left(\begin{array}{cc}
0 &  1\\
1 & 0\end{array}\right),\qquad A_{12}^{\pm, 1}=A_{21}^{\pm,
1}=A_{22}^{\pm, 1}=0
\end{equation}
in a neighborhood of $\{x_1=0\}$.

Thus, from (\ref{f3.10}) and (\ref{f4.2}), we find that the vector
function
$$
Y^\pm=(P^\pm)^{-1}X^\pm
$$
satisfies the following problem:
%\stepcounter{line}
\begin{equation}\label{f4.4}
\begin{cases}{\mathcal L}(U^\pm, \Psi^\pm)Y^\pm+{\mathcal E}^\pm(U^\pm,
\Psi^\pm)Y^\pm={\cal F}^\pm \qquad\mbox{in}\,\, \{x_1>0\},\cr
[Y_1]=[Y_2]=0 \qquad{\rm on}\,\, \{x_1=0\},\cr Y^\pm|_{t\le 0}=0,
\end{cases}
\end{equation}
where ${\cal F}^\pm=(P^\pm)^\top F^\pm$,
$$
{\cal L}(U^\pm, \Psi^\pm) ={\cal A}_0^\pm \partial_t
+\sum_{j=1}^3{\cal A}_j^\pm\partial_{x_j}
$$
with ${\cal A}_j^\pm=(P^\pm)^\top A_j^\pm P^\pm, 0\le j\le 3$, and $
{\cal E}^\pm=(P^\pm)^\top E P^\pm+ (P^\pm)^\top \sum_{j=0}^3
A_j^\pm \partial_{x_j}P^\pm.$

\medskip
{\it $L^2$-estimate}: With setup \eqref{f4.4}, we can apply the
Lax-Friedrichs theory to obtain the following $L^2-$estimate:
\begin{equation}\label{f4.5}
\max_{0\le t\le T}\|Y^\pm(t)\|_{0, \mu}^2+
\mu\|Y^\pm\|_{0, \mu, T}^2
\le \frac{C_0}{\mu}\|{\cal F}^\pm\|_{0, \mu, T}^2
\end{equation}
when $\mu\ge \mu_0>0$, where
$$
C_0=C_0(\|coef\|_{L^\infty}, \|\nabla\cdot coef\|_{L^\infty}).
$$

\medskip
{\it Estimates on the tangential derivatives}. From (\ref{f4.4}), we
find that $M^\alpha Y^\pm, |\alpha|\le s,$ satisfy
%\stepcounter{line}
\begin{equation}\label{f4.6}
\begin{cases}
{\cal L}M^\alpha Y^\pm+{\cal E}^\pm M^\alpha Y^\pm=M^\alpha {\cal
F}^\pm+[ {\cal L}+{\cal E}^\pm, M^\alpha]Y^\pm \qquad\mbox{in}\,\,
\{x_1>0\},\cr [M^\alpha Y_1]=[M^\alpha Y_2]=0\qquad{\rm on}\,\,
\{x_1=0\},\cr M^\alpha Y^\pm|_{t\le 0}=0,
\end{cases}
\end{equation}
which, by the Lax-Friedrichs theory again, yields
%\stepcounter{line}
\begin{equation}\label{f4.7}
\begin{aligned}
&\max_{0\le t\le T}|Y^\pm(t)|_{s, \mu}^2+ \mu|Y^\pm|_{s, \mu, T}^2\\
&\le \frac{C_0}{\mu}\big(|{\cal F}^\pm|_{s, \mu,
T}^2+(\|coef\|_T^\ast)^2 |Y^\pm|_{s, \mu, T}^2+\|\nabla
Y^\pm\|^2_{L^\infty}|\dot{coef}|_{s, \mu,
T}^2+|\partial_{x_1}Y^\pm_{I}|_{s-1, \mu, T}^2\big),
\end{aligned}
\end{equation}
where $Y^\pm_I=(Y_1^\pm, Y_2^\pm)$, the norms $|\cdot|_{s,\mu}$ and
$|\cdot|_{s,\mu, T}$ are the special cases of the norms
$\|\cdot\|_{s,\mu}$ and $\|\cdot\|_{s,\mu, T}$ respectively without
the normal derivatives $\partial_{x_1}^k$,
and
$$\|u\|_{T}^\ast=\sum_{|\alpha|\le 2}\|M^\alpha u\|_{L^\infty(\Omega_T)}+
\sum_{|\alpha|\le 1}\|\partial_{x_1}M^\alpha u\|_{L^\infty(\Omega_T)}.$$

\medskip
{\it Estimates on the normal derivatives}. Set ${\cal
A}_j=diag[{\cal A}_j^+, {\cal A}_j^-], 0\le j\le 3,$ with
$${\cal A}_j^\pm=\left(\begin{array}{cc}
A_{11}^{\pm, j} &  A_{12}^{\pm, j}\\
A_{21}^{\pm, j} &A_{22}^{\pm, j}\end{array}\right).
$$
{}From (\ref{f4.3})--(\ref{f4.4}), we find that $Y_I^\pm=(Y_1^\pm,
Y_2^\pm)$ and $Y_{I\!\!I}^\pm=(Y_3^\pm, \ldots, Y_8^\pm)$ satisfy
the following equations:
%\stepcounter{line}
\begin{equation}\label{f4.8}
A_{11}^{\pm, 1}\partial_{x_1}Y_I^\pm= coef\cdot M Y^\pm+ coef\cdot
Y^\pm +{\cal F}_I^\pm,
\end{equation}
and
\begin{equation}\label{f4.9}
\begin{cases} \sum_{j=0}^3A_{22}^{\pm,
j}\partial_{x_j}Y_{I\!\!I}^\pm+C_{22}^{\pm}Y_{I\!\!I}^\pm ={\cal
F}_{I\!\!I}^\pm-\sum_{j=0}^3A_{21}^{\pm, j}\partial_{x_j}Y_{I}^\pm-
C_{21}^{\pm}Y_{I}^\pm,\cr Y_{I\!\!I}^\pm|_{t\le 0}=0,
\end{cases}
\end{equation}
where the linear operator
$\sum_{j=0}^3A_{22}^{\pm,j}\partial_{x_j}+C_{22}^{\pm}$ is symmetric
hyperbolic and tangential to $\{x_1=0\}$.

Thus, from (\ref{f4.8}), we obtain
\begin{equation}\label{f4.10}
|\partial_{x_1}Y_I^\pm(t)|_{s-1, \mu} \le |{\cal F}_I^\pm(t)|_{s-1,
\mu}+ \|coef\|_{Lip(T)}|Y^\pm(t)|_{s, \mu}
+\|Y^\pm\|_{Lip(T)}|\dot{coef}(t)|_{s, \mu},
\end{equation}
with $\|\cdot\|_{Lip(T)}$ denoting the usual norm in the space of
Lipschitz continuous functions in $\Omega_T$ and, from (\ref{f4.9}),
we deduce the following estimate on $\partial_{x_1}Y_{I\!\!I}^\pm$
by using the classical hyperbolic theory:
\begin{equation}\label{f4.11}
\begin{aligned}
&\max_{0\le t\le T}|\partial_{x_1}Y_{I\!\!I}^\pm(t)|_{s-2, \mu}^2
+\mu |\partial_{x_1}Y_{I\!\!I}^\pm|_{s-2, \mu, T}^2\\
& \le \frac{C_0}{\mu}\Big(|\partial_{x_1}{\cal
F}^\pm_{I\!\!I}|_{s-2, \mu, T}^2+ \|coef\|_{Lip(T)}^2(|Y^\pm_I|_{s,
\mu, T}^2
+|\partial_{x_1}Y^\pm_{I}|_{s-1, \mu, T}^2)\\
&\quad\qquad +(\|coef\|_{T}^\ast|Y^\pm_{I\!\!I}|_{s, \mu, T}^N)^2+
(\|Y^\pm\|^\ast_T|\dot{coef}|^N_{s, \mu, T})^2\Big),
\end{aligned}
\end{equation}
where
$$
|u|^N_{s,\mu,T}=|u|_{s,\mu,T}+|u_{x_1}|_{s-2,\mu,T}.
$$

Combining estimates (\ref{f4.7}) with (\ref{f4.10})--(\ref{f4.11}),
applying an inductive argument on the order of normal derivatives
$\partial_{x_1}^k$, and noting that all the coefficient functions of
$\cal L$ and ${\cal E}^\pm$ in problem (\ref{f4.4}) depending on
$(U^\pm, \nabla_{t,x}\Psi^\pm)$ and $(U^\pm, \nabla_{t,x}U^\pm,
\{\partial_{x_1}^{\alpha_1}\nabla_{t,x_2,
x_3}^{\alpha_2}\Psi^\pm\}_{\alpha_1\le 1, |\alpha_2|\le 2,
\alpha_1+|\alpha_2|\le 2})$, respectively, we can conclude estimate
(\ref{f4.1}) for the case of homogeneous boundary conditions by
choosing $s_0> 17/2$.

\medskip
{\it Step 2}. With Step 1, for the case of nonhomogeneous boundary
conditions in (\ref{f3.10}) with $(h_1^\pm, h_2^\pm, h_3)\in
H^{s+1}([0,T]\times \R^2)$, it suffices to study (\ref{f3.10}) when
$F^\pm=0$.

By Lemma 4.1, there exists $X_0^\pm\in B^{s+2}(\Omega_T)$ satisfying
%\stepcounter{line}
\begin{equation}\label{f4.12}
[X_{0,2}-\lambda X_{0,5}]=-[h_1+\lambda h_2], \quad [X_{0,1}]=h_3
\qquad {\rm on}\,\, \{x_1=0\}.
\end{equation}

Thus, from  (\ref{f3.10}), we know that $Y^\pm=X^\pm-X^\pm_0$
satisfies the following problem:
%\stepcounter{line}
\begin{equation}\label{f4.13}
\begin{cases}
L(U^\pm, \Psi^\pm)Y^\pm+E(U^\pm, \Psi^\pm)Y^\pm=f^\pm\qquad
\mbox{in}\,\, \{x_1>0\},\cr
[Y_1]=[Y_2-\lambda Y_5]=0 \qquad
\mbox{on}\,\, \{x_1=0\},\cr
 Y^\pm|_{t\le 0}=0, \end{cases}
\end{equation}
with $f^\pm=-L(U^\pm, \Psi^\pm)X_0^\pm-E(U^\pm, \Psi^\pm)X_0^\pm\in
B^s(\Omega_T)$.

Employing estimate (\ref{f4.1}) in the case $h=0$ for problem
(\ref{f4.13}) established in Step 1 yields an estimate on $Y^\pm$,
which implies the estimates of $X^\pm=Y^\pm+X_0^\pm$.
Combining this estimate with that for the case of homogeneous
boundary conditions yields (\ref{f4.1}) for the general case.

\vspace{.1in} By choosing $\mu\gg 1$ and $T$ so that $\mu T\le 1$ in
(\ref{f4.1}), we can directly conclude

\vspace{.1in} {\bf Theorem 4.2.} {\it For any fixed $s_0>{17}/{2}$,
there exists a constant $C_0>0$ depending only on
$\|\dot{coef}\|_{s_0,T}$ such that, for any $s\ge s_0$, the estimate
\begin{equation}\label{f4.14}
\|X^\pm\|_{s, T}^2 \le C_0\big(\|{F}^\pm\|_{s, T}^2
+\|{h}\|_{H^{s+1}(b\Omega_T)}^2 +\|\dot{coef}\|_{s, T}^2
(\|{F}^\pm\|_{s_0, T}^2
+\|{h}\|_{H^{s_0+1}(b\Omega_T)}^2)\big)
\end{equation}
holds.}

\section{Nonlinear Stability I: Construction of The Zero-th Order Approximate
Solutions}

With the linear estimates in Section 4, we now establish the local
existence of current-vortex sheets of problem
(\ref{f2.15})--(\ref{f2.18}) under the compatibility conditions on
the initial data. From now on, we focus on the initial data
$(U^\pm_0, \psi_0)$ that is a small perturbation of a planar
current-vortex sheet $(\bar{U}^\pm, \bar{\psi})$ for constant states
$\bar{U}^\pm$ and $\bar{\psi}=0$ satisfying
\eqref{f2.3}--\eqref{f2.4}.

\medskip
We start with {\it the compatibility conditions on $\{t=x_1=0\}$}.
For fixed $k\in \N$, it is easy to formulate the $j^{th}$--order
compatibility conditions, $0\le j\le k$, for the initial data
$(U^\pm_0, \psi_0)$ of problem \eqref{f2.15}--(\ref{f2.18}), from
which we determine the data:
%\stepcounter{nonl}
\begin{equation}\label{f5.1}
\{\partial^j_t U^\pm|_{t=x_1=0}\}_{0\le j\le k}\quad {\rm and}\quad
\{\partial^j_t \psi|_{t=0}\}_{0\le j\le k+1}.
\end{equation}

For any fixed integer $s>9/2$ and given data $\psi_0\in
H^{s-1}(\R^2)$ and $U^\pm_0$ with
$$
\dot{U}_0^\pm={U}_0^\pm-\bar{U}^\pm\in B^s(\R_+^3),
$$
by using Lemma 4.1, we can extend $\psi_0$ to be
$\dot{\Psi}^\pm_0\in B^{s}(\R_+^3)$ satisfying
$$U_{0,5}^\pm-\partial_{x_2}\dot{\Psi}_0^\pm U_{0,6}^\pm
-\partial_{x_3}\dot{\Psi}_0^\pm U_{0,7}^\pm=0$$ which is possible by
using assumption (\ref{f2.3a}). Let $\Psi_0^\pm=\dot{\Psi}^\pm_0\pm
x_1$ be the initial data for solving $\Psi^\pm$ from the eikonal
equations (\ref{f2.15}).

Let $U_j^\pm(x)=\partial^j_t U^\pm|_{t=0}$ and
$\Psi_j^\pm=\partial^j_t \Psi^\pm|_{t=0}$ for any $j\ge 0$.

For all $k\le [s/2]$, from problem (\ref{f2.15})--(\ref{f2.18}), we
determine
\begin{equation}\label{f5.2}
\{
U_j^\pm\}_{1\le j\le k}\quad {\rm and}\quad \{\Psi_j^\pm\}_{0\le j\le k+1}
\end{equation}
in terms of $(U^\pm_0, \Psi_0^\pm)$ and
\begin{equation}\label{f5.3}
\begin{cases}
U_j^\pm\in B^{s-2j}(\R^3_+), \quad 1\le j\le k,\\
\Psi^\pm_1\in B^{s-1}(\R^3_+), \quad \Psi^\pm_j\in
B^{s-2(j-1)}(\R^3_+), \,\, 2\le j\le k+1.
\end{cases}
\end{equation}
Furthermore, we have the estimate
\begin{equation}\label{f5.4}
\begin{aligned}
&\sum_{j=1}^{[s/2]}\|U_j^\pm\|_{B^{s-2j}(\R_+^3)}
+\sum_{j=2}^{[s/2]+1}\|\Psi_j^\pm\|_{B^{s-2j+2}(\R_+^3)}+
\|\Psi_1^\pm\|_{B^{s-1}(\R_+^3)}\\
&\le C_0
(\|\dot{U}_0^\pm\|_{B^s(\R^3_+)}+\|\psi_0\|_{H^{s-1}(\R^2)})
\end{aligned}
\end{equation}
with $C_0$ depending only on $s$,
$\|\dot{U}_0^\pm\|_{W^{1,\infty}(\R_+^3)}$, and
$\|\psi_0\|_{W^{1,\infty}(\R^2)}$.

\bigskip Set $s_1=[s/2]+1$. We now construct the {\it zero-th order
approximate solutions $(U_a^\pm, \Psi_a^\pm)$}. First, we construct
\begin{equation}\label{f5.5}
\dot{\Psi}_a^\pm\in B^{s_1+1}(\R_+\times \R_+^3)
\end{equation}
satisfying
$\partial_t^j\dot{\Psi}_a^\pm|_{t=0}=\dot{\Psi}_j^\pm, 0\le j\le
s_1$, and $\dot{\Psi}_a^+|_{x_1=0}=\dot{\Psi}_a^-|_{x_1=0}$,
and
\begin{equation}\label{f5.7}
\dot{U}_{a, l}^\pm\in B^{s_1}(\R_+\times \R_+^3), \quad l\neq 2,
\end{equation}
satisfying
$\partial_t^j\dot{U}_{a,l}^\pm|_{t=0}=\dot{U}_{j,l}^\pm, 0\le j\le
s_1-1$.

Moreover, we have
\begin{equation}\label{f5.9}
\begin{cases}
U^\pm_{a,5}-\partial_{x_2}\Psi_a^\pm U^\pm_{a,6}-
\partial_{x_3}\Psi_a^\pm U^\pm_{a,7}=0\qquad \mbox{in}\,\, \{x_1>0\},\cr
[U_{a,1}+\frac 12 |U_{a,H}|^2]=0 \qquad \mbox{on}\,\,
\{x_1=0\},\end{cases}
\end{equation}
where $U_a^\pm=\bar{U}^\pm+\dot{U}^\pm_a$,
$\Psi_a^\pm=\dot{\Psi}^\pm_a\pm x_1$,  and $U_{a, H}=(U_{a, 5},
U_{a, 6}, U_{a,7})^\top$.

Finally, we construct $U_{a,2}^\pm\in B^{s_1}(\R_+\times \R_+^3)$ by
requiring
\begin{equation}\label{f5.10}
\partial_t\Psi_a^\pm-U^\pm_{a,2}+
\partial_{x_2}\Psi_a^\pm U^\pm_{a,3}+
\partial_{x_3}\Psi_a^\pm U^\pm_{a,4}=0 \qquad
{\rm in}~\{x_1>0\}.
\end{equation}

Therefore, we conclude

\vspace{.1in} {\bf Lemma 5.1.} {\it For any fixed integer $s>9/2$,
assume that the initial data
$(\dot{U}^\pm_0,\psi_0)=({U}_0^\pm-\bar{U}^\pm, \psi_0)\in
B^s(\R_+^3)\times H^{s-1}(\R^2)$ is bounded in the norms. Then there
exists $(U_a^\pm, \Psi_a^\pm)$ such that
$\dot{U}_a^\pm={U}_a^\pm-\bar{U}^\pm\in B^{[s/2]+1}(\R_+\times
\R_+^3)$, $\dot{\Psi}_a^\pm={\Psi}_a^\pm\mp x_1\in
B^{[s/2]+2}(\R_+\times \R_+^3)$,
\begin{equation}\label{f5.11}
\partial_t^j(L(U_a^\pm, \Psi_a^\pm)U_a^\pm))|_{t=0}=0 \quad\mbox{for}\,\, 0\le j\le
[s/2]-1,
\end{equation}
and {\rm (\ref{f5.9})--(\ref{f5.10})} hold.}

\vspace{.1in} With this, we set $\psi_a=\Psi_a^\pm|_{x_1=0}$. From
the compatibility conditions of problem \eqref{f2.15}--(\ref{f2.18})
and Lemma 5.1, we have
\begin{equation}\label{f5.12}
\partial_t^j(B(U_a^+, U_a^-, \psi_a))|_{t=0}=0 \qquad\mbox{for}\,\, 0\le j\le [s/2].
\end{equation}

\vspace{.15in} Set
\begin{equation}\label{f3.13}
 V^\pm=U^\pm-U^\pm_a,\qquad \Phi^\pm=\Psi^\pm-\Psi_a^\pm.
\end{equation}
Then it follows from Lemma 5.1 and (\ref{f5.12}) that problem
(\ref{f2.15})--(\ref{f2.18}) is equivalent to the following {\bf
fixed initial-boundary value problem for $(V^\pm, \Phi^\pm)$}:
\begin{equation}\label{f5.14}
\begin{cases}
{\cal L}(V^\pm, \Phi^\pm)V^\pm=f_a^\pm \qquad {\rm in}~ \{t>0,
~x_1>0\},\cr {\cal E}(V^\pm, \Phi^\pm)=0 \qquad {\rm in}~ \{t>0,
~x_1>0\},\cr \Phi^+|_{x_1=0}=\Phi^-|_{x_1=0}=\phi,\cr
{\cal
B}(V^+,V^-, \phi)=0\qquad {\rm on}~ \{x_1=0\},\cr V^\pm|_{t\le 0}=0,
\quad \Phi^\pm|_{t\le 0}=0,
\end{cases}
\end{equation}
where $f_a^\pm=-L(U^\pm_a, \Psi_a^\pm)U_a^\pm$,
\begin{eqnarray*}
&& {\cal L}(V^\pm, \Phi^\pm)V^\pm=L(U_a^\pm+V^\pm,
\Psi_a^\pm+\Phi^\pm)(U_a^\pm+V^\pm)- L(U_a^\pm,
\Psi_a^\pm)U_a^\pm,\\
&&{\cal E}(V^\pm, \Phi^\pm)
=\partial_t\Phi^\pm-V_2^\pm+\partial_{x_2}(\Psi_a^\pm+\Phi^\pm)V_3^\pm
+\partial_{x_3}(\Psi_a^\pm+\Phi^\pm)V_4^\pm
+U_{a,3}^\pm\partial_{x_2}\Phi^\pm+U_{a,4}^\pm\partial_{x_3}\Phi^\pm,
\end{eqnarray*}
and
$$
{\cal B}(V^+,V^-, \phi)=B(U_a^++V^+,U_a^-+V^-, \psi_a+\phi).
$$

\section{Nonlinear Stability II: Iteration Scheme}

{}From Theorem 4.1, we observe that the high-order energy estimates
of solutions for the linearized problem (\ref{f3.10})--\eqref{f3.11}
in terms of the nonhomogeneous terms $F^\pm$ and variable
coefficients $\dot{coef}$ keep the same order.
Based on these estimates and the structure of the nonlinear system,
we now develop a suitable iteration scheme of
Nash-Moser-H\"{o}rmander type (cf. \cite{hor}) for our {\bf
nonlinear problem (\ref{f5.14})}.

To do this, we first recall a standard family of smoothing operators
(cf. \cite{ali, cou2}):
%\stepcounter{nonl}
\begin{equation}\label{f6.1}
\{S_\theta\}_{\theta>0}:
~ B_\mu^0(\Omega_T)\longrightarrow \cap_{s\ge 0}B^s_\mu(\Omega_T)
\end{equation}
satisfying
\begin{equation}\label{f6.2}
\begin{cases}
\|S_\theta u\|_{s, T}\le C\theta^{(s-\alpha)_+}\|u\|_{\alpha, T}
\quad {\rm for ~all}~s, \alpha\ge 0,\cr \|S_\theta u-u\|_{s, T}\le
C\theta^{s-\alpha}\|u\|_{\alpha, T} \quad {\rm for ~all}~ s\in [0,
\alpha],\cr \|\frac{d}{d\theta}S_\theta u\|_{s, T}\le
C\theta^{s-\alpha-1}\|u\|_{\alpha, T} \quad {\rm for ~all}~s,
\alpha\ge 0, \end{cases}
\end{equation}
and
\begin{equation}\label{f6.3}
\|(S_\theta u_+-S_\theta u_-)|_{x_1=0}\|_{H^s(b\Omega_T)} \le
C\theta^{(s+1-\alpha)_+}\|(u_+-u_-)|_{x_1=0}\|_{\alpha, T} \qquad
{\rm for ~all}~s, \alpha\ge 0.
\end{equation}

Similarly, one has a family of smoothing operators, still denoted by
$\{S_\theta\}_{\theta>0}$ acting on $H^s(b\Omega_T)$, satisfying
(\ref{f6.2}) as well for the norms of $H^s(b\Omega_T)$ (cf.
\cite{ali, cou2}).

Now we construct the iteration scheme for solving the nonlinear
problem \eqref{f5.14} in $\R_+\times\R^3$.

\medskip
{\sc The Iteration Scheme:} Let $V^{\pm, 0}=\Phi^{\pm, 0}=0$. Assume
that $(V^{\pm, k}, \Phi^{\pm, k})$ have been known for $k=0, \ldots,
n$, and satisfy
\begin{equation}\label{f6.4}
\begin{cases}(V^{\pm, k}, \Phi^{\pm, k})=0\qquad {\rm in}~\{t\le 0\},\cr
\Phi^{+, k}|_{x_1=0}=\Phi^{-, k}|_{x_1=0}=\phi^k\quad {\rm
on}~\{x_1= 0\}.
\end{cases}
\end{equation}
Denote the $(n+1)^{th}$ approximate solutions to (\ref{f5.14}) in
$\R_+\times\R^3$ by
\begin{equation}\label{f6.5}
V^{\pm, n+1}=V^{\pm, n}+\delta V^{\pm, n},\quad \Phi^{\pm,
n+1}=\Phi^{\pm, n}+\delta\Phi^{\pm, n}, \quad
\phi^{n+1}=\phi^n+\delta\phi^n.
\end{equation}

\medskip
Let $\theta_0\ge 1$ and $\theta_n=\sqrt{\theta_0^2+n}$ for any $n\ge
1$. Let $S_{\theta_n}$ be the associated smoothing operators defined
as above. We now determine the problem of the increments $(\delta
V^{\pm, n}, \delta\Phi^{\pm, n})$ as follows:

\medskip
{\it Construction of $\delta V^{\pm,n}$}: First, it is easy to see
\begin{equation}\label{f6.6}
{\cal L}(V^{\pm, n+1}, \Phi^{\pm, n+1})V^{\pm,n+1}-{\cal L}(V^{\pm,
n},\Phi^{\pm, n})V^{\pm, n} =L'_{e, (U_a^\pm+V^{\pm, n+\frac 12},
\Psi_a^\pm+S_{\theta_n}\Phi^{\pm, n})} \delta \dot{V}^{\pm, n}
+e_{\pm, n},
\end{equation}
where $V^{\pm, n+\frac 12}$ are the modified states of
$S_{\theta_n}V^{\pm, n}$ defined in (\ref{f6.24})--(\ref{f6.25})
later, which guarantees that the boundary $\{x_1=0\}$ is the uniform
characteristic of constant multiplicity at each iteration step,
\begin{equation}\label{f6.7}
\begin{aligned}
&L'_{e, (U_a^\pm+V^{\pm, n+\frac 12},
\Psi_a^\pm+S_{\theta_n}\Phi^{\pm, n})} \delta \dot{V}^{\pm,
n}\\
&=L(U_a^\pm+V^{\pm, n+\frac 12}, \Psi_a^\pm+S_{\theta_n}\Phi^{\pm,
n})\delta \dot{V}^{\pm, n} + E(U_a^\pm+V^{\pm, n+\frac 12},
\Psi_a^\pm+S_{\theta_n}\Phi^{\pm, n})\delta \dot{V}^{\pm, n}
\end{aligned}
\end{equation}
are the effective linear operators,
\begin{equation}\label{f6.8}
\delta \dot{V}^{\pm, n}=\delta {V^{\pm, n}}
-\delta \Phi^{\pm, n}\frac{\partial_{x_1}(U_a^\pm+V^{\pm, n+\frac 12})}
{\partial_{x_1}(\Psi_a^\pm+S_{\theta_n}\Phi^{\pm, n})}
\end{equation}
are the good unknowns, and
\begin{equation}\label{f6.9}
e_{\pm, n}=\sum_{j=1}^4e_{\pm, n}^{(j)}
\end{equation}
are the total errors with the first error resulting from the Newton
iteration scheme:
\begin{equation}\label{f6.10}
\begin{aligned}
e_{\pm, n}^{(1)}=&L(U_a^\pm+V^{\pm, n+1}, \Psi_a^\pm+\Phi^{\pm,
n+1})(U_a^\pm+V^{\pm, n+1})-L(U_a^\pm+V^{\pm, n},
\Psi_a^\pm+\Phi^{\pm, n})(U_a^\pm+V^{\pm,
n})\\
& -L'_{(U_a^\pm+V^{\pm, n}, \Psi_a^\pm+\Phi^{\pm, n})}(\delta
{V^{\pm, n}}, \delta {\Phi^{\pm, n}}),
\end{aligned}
\end{equation}
the second error resulting from the substitution:
\begin{equation}\label{f6.11}
e_{\pm, n}^{(2)} =L'_{(U_a^\pm+V^{\pm, n},\Psi_a^\pm+\Phi^{\pm,
n})}(\delta {V^{\pm, n}}, \delta {\Phi^{\pm, n}})
-L'_{(U_a^\pm+S_{\theta_n}V^{\pm, n},
\Psi_a^\pm+S_{\theta_n}\Phi^{\pm, n})}(\delta {V^{\pm, n}}, \delta
{\Phi^{\pm, n}}),
\end{equation}
the third error resulting from the second substitution:
\begin{equation}\label{f6.12}
e_{\pm, n}^{(3)}= L'_{(U_a^\pm+S_{\theta_n}V^{\pm, n}, \Psi_a^\pm
+S_{\theta_n}\Phi^{\pm, n})}(\delta {V^{\pm, n}}, \delta {\Phi^{\pm,
n}}) -L'_{(U_a^\pm+V^{\pm, n+\frac 12},
\Psi_a^\pm+S_{\theta_n}\Phi^{\pm, n})}(\delta {V^{\pm, n}}, \delta
{\Phi^{\pm, n}}),
\end{equation}
and the remaining error:
\begin{equation}\label{f6.13}
e_{\pm, n}^{(4)}= \frac{\delta \Phi^{\pm,
n}}{\partial_{x_1}\big(\Psi_a^\pm+S_{\theta_n}\Phi^{\pm, n})}
\partial_{x_1}(L(U_a^\pm+V^{\pm, n+\frac 12},
\Psi_a^\pm+S_{\theta_n}\Phi^{\pm, n})(U_a^\pm+V^{\pm, n+\frac
12})\big).
\end{equation}

Similarly, we have
\begin{equation}\label{f6.14}
{\cal B}(V^{+, n+1}, V^{-, n+1}, \phi^{n+1})-{\cal B}(V^{+, n},
V^{-, n}, \phi^{n})=B'_{(U_a^\pm+V^{\pm, n+\frac 12},
\psi_a+S_{\theta_n}\phi^{n})}(\delta \dot{V}^{+, n}, \delta
\dot{V}^{-, n},\delta \phi^{n})+\tilde{e}_n,
\end{equation}
where
%\stepcounter{nonl}
\begin{equation}\label{f6.15}
\tilde{e}_{n}=\sum_{j=1}^4\tilde{e}_{n}^{(j)}
\end{equation}
with the first error resulting from the Newton iteration scheme:
%\stepcounter{nonl}
\begin{equation}\label{f6.16}
\begin{aligned} \tilde{e}_{ n}^{(1)}=&B(U_a^++V^{+, n+1}, U_a^-+V^{-,
n+1},\psi_a+\phi^{n+1})- B(U_a^++V^{+, n}, U_a^-+V^{-,
n},\psi_a+\phi^{n})\\
& -B'_{(U_a^\pm+V^{\pm, n}, \psi_a+\phi^{n})}(\delta {V^{+,
n}},\delta {V^{-, n}}, \delta {\phi^{n}}),
\end{aligned}
\end{equation}
the second error resulting from the substitution:
\begin{equation}\label{f6.17}
\tilde{e}_{n}^{(2)} =B'_{(U_a^\pm+V^{\pm, n}, \psi_a+\phi^{ n})}
(\delta {V^{+, n}},\delta {V^{-, n}},
 \delta {\phi^{n}})
-B'_{(U_a^\pm+S_{\theta_n}V^{\pm, n},
\psi_a+S_{\theta_n}\phi^{n})}(\delta {V^{+, n}},\delta {V^{-, n}},
 \delta {\phi^{n}}),
\end{equation}
the third errors resulting from the second substitution:
\begin{equation}\label{f6.18}
\tilde{e}_{n}^{(3)}= B'_{(U_a^\pm+S_{\theta_n}V^{\pm, n},
\psi_a+S_{\theta_n}\phi^{n})}(\delta {V^{+, n}},\delta {V^{-, n}},
 \delta {\phi^{n}})
-B'_{(U_a^\pm+V^{\pm, n+\frac 12}, \psi_a +S_{\theta_n}\phi^{n})}
(\delta {V^{+, n}},\delta {V^{-, n}},
 \delta {\phi^{n}}),
\end{equation}
and the remaining error:
\begin{equation}\label{f6.19}
\tilde{e}_{n}^{(4)}= B'_{(U_a^\pm+V^{\pm, n+\frac 12},
\psi_a+S_{\theta_n}\phi^{n})}(\delta {V^{+, n}},\delta {V^{-, n}},
 \delta {\phi^{n}})-B'_{(U_a^\pm+V^{\pm, n+\frac 12}, \psi_a+S_{\theta_n}\phi^{n})}
(\delta {\dot{V}^{+, n}},\delta {\dot{V}^{-, n}},
 \delta {\phi^{n}}).
\end{equation}

\medskip
Using (\ref{f6.6}) and (\ref{f6.14}) and noting that
$$
%\begin{cases}
{\cal L}(V^{\pm, 0}, \Phi^{\pm, 0})V^{\pm, 0}=0, \qquad {\cal
B}(V^{+, 0},V^{-, 0}, \phi^{0})=0,
%\end{cases}
$$
we obtain
%\stepcounter{nonl}
\begin{equation}\label{f6.20}
\begin{cases} {\cal L}(V^{\pm, n+1}, \Phi^{\pm, n+1})V^{\pm,
n+1}=\sum_{j=0}^n\big( L'_{e,(U_a^\pm+V^{\pm, j+\frac 12},
\Psi_a^\pm+S_{\theta_j}\Phi^{\pm,j})}\delta \dot{V}^{\pm, j}+e_{\pm,
j}\big), \\
{\cal B}(V^{+, n+1},V^{-, n+1}, \phi^{n+1})=\sum_{j=0}^n\big(
B'_{(U_a^\pm+V^{\pm, j+\frac 12},
\psi_a+S_{\theta_j}\phi^{j})}(\delta \dot{V}^{+, j}, \delta
\dot{V}^{-, j}, \delta\phi^j)+\tilde{e}_{j}\big).
\end{cases}
\end{equation}

Observe that, if the limit of $(V^{\pm, n}, \Phi^{\pm, n})$ exists
which is expected to be a solution to problem (\ref{f5.14}), then
the left-hand side in the first equation of (\ref{f6.20}) should
converge to $f_a^\pm$, and the left-hand side in the second one of
(\ref{f6.20}) goes to zero when $n\to \infty$. Thus, from
(\ref{f6.20}), it suffices to study the following problem:
% \stepcounter{nonl}
\begin{equation}\label{f6.21}
\begin{cases}
L'_{e,(U_a^\pm+V^{\pm, n+\frac 12},
\Psi_a^\pm+S_{\theta_n}\Phi^{\pm,n})}\delta \dot{V}^{\pm, n}
=f_n^\pm \qquad {\rm in}~\Omega_T,\cr B'_{(U_a^\pm+V^{\pm, n+\frac
12}, \psi_a+S_{\theta_n}\phi^{n})}(\delta \dot{V}^{+, n}, \delta
\dot{V}^{-, n}, \delta\phi^n)=g_n \qquad {\rm on}~b\Omega_T,\cr
\delta \dot{V}^{\pm, n}|_{t\le 0}=0, \end{cases}
\end{equation}
where
$f_n^\pm$ and $g_n$ are defined by
%\stepcounter{nonl}
\begin{equation}\label{f6.22}
%\begin{cases}
 \sum_{j=0}^nf_j^\pm+S_{\theta_n}(\sum_{j=0}^{n-1}e_{\pm,
j})=S_{\theta_n}f_a^\pm, \qquad
\sum_{j=0}^ng_j+S_{\theta_n}(\sum_{j=0}^{n-1}\tilde{e}_{j})=0,
%\end{cases}
\end{equation}
by induction on $n$, with
$f_0^\pm=S_{\theta_0}f_a^\pm$ and $g_0=0$.

\medskip
We now define the modified state $V^{\pm, n+\frac 12}$ to guarantee
that the boundary $\{x_1=0\}$ is the uniform characteristic of
constant multiplicity at each iteration step (\ref{f6.21}). To
achieve this, we require
%\stepcounter{nonl}
\begin{equation}\label{f6.23}
({\cal B}(V^{+, n+\frac 12},V^{-, n+\frac 12},
S_{\theta_n}\phi^{n}))_i^\pm|_{x_1=0}=0 \qquad \mbox{for}\,\, i=1,2,
\,\,\mbox{and all}\,\, n\in \N,
\end{equation}
which leads to define
\begin{equation}\label{f6.24}
V_j^{+, n+\frac 12}=S_{\theta_n}V^{\pm, n}_j\qquad\mbox{for}\,\,
j\in \{3,4,6,7\},
\end{equation}
 and
%\stepcounter{nonl}
\begin{equation}\label{f6.25}
\begin{cases} V_2^{\pm, n+\frac 12}=\partial_t(S_{\theta_n}\Phi^{\pm,
n})+
\partial_{x_2}(\Psi^{\pm}_a+S_{\theta_n}\Phi^{\pm, n})V_3^{\pm, n+\frac 12}
+\partial_{x_3}(\Psi^{\pm}_a+S_{\theta_n}\Phi^{\pm, n})V_4^{\pm,
n+\frac 12}\\ \hspace{.7in}
+U^\pm_{a,3}\partial_{x_2}(S_{\theta_n}\Phi^{\pm,
n})+U^\pm_{a,4}\partial_{x_3}(S_{\theta_n}\Phi^{\pm, n}), \\
V_5^{\pm, n+\frac
12}=\partial_{x_2}(\Psi^{\pm}_a+S_{\theta_n}\Phi^{\pm, n})V_6^{\pm,
n+\frac 12} +\partial_{x_3}(\Psi^{\pm}_a+S_{\theta_n}\Phi^{\pm,
n})V_7^{\pm, n+\frac 12}\\ \hspace{.7in}
+U^\pm_{a,6}\partial_{x_2}(S_{\theta_n}\Phi^{\pm,
n})+U^\pm_{a,7}\partial_{x_3}(S_{\theta_n}\Phi^{\pm, n}).
\end{cases}
\end{equation}

\bigskip {\it Construction of $\delta\Phi^{\pm,n}$ satisfying
$\delta\Phi^{\pm,n}|_{x_1=0}=\delta\phi^{n}$}. Clearly, we have
%\stepcounter{nonl}
\begin{equation}\label{f6.26}
{\cal E}(V^{\pm, n+1}, \Phi^{\pm, n+1})-{\cal E}(V^{\pm,
n},\Phi^{\pm, n})= {\cal E}'_{(V^{\pm, n+\frac 12},
S_{\theta_n}\Phi^{\pm,n})}(\delta \dot{V}^{\pm, n}, \delta
\Phi^{\pm,n})+\bar{e}_{\pm,n},
\end{equation}
where
%\stepcounter{nonl}
\begin{equation}\label{f6.27}
\begin{aligned}
{\cal E}'_{(V^{\pm}, \Phi^{\pm})}(W^{\pm}, \Theta^{\pm})=&
\partial_t\Theta^{\pm}-W_2^{\pm}+\partial_{x_2}(\Psi^\pm_a+\Phi^{\pm})
W_3^{\pm} +\partial_{x_3}(\Psi^\pm_a+\Phi^{\pm})W_4^{\pm}\cr
&+(U_{a,3}^\pm+V_3^{\pm})\partial_{x_2}\Theta^{\pm}
+(U_{a,4}^\pm+V_4^{\pm})\partial_{x_3}\Theta^{\pm}
\end{aligned}
\end{equation}
is the linearized operator of $\cal E$ and
%\stepcounter{nonl}
\begin{equation}\label{f6.28}
\bar{e}_{\pm,n}=\sum_{j=1}^4\bar{e}_{\pm, n}^{(j)}
\end{equation}
with the first error resulting from the Newton iteration scheme:
%\stepcounter{nonl}
\begin{equation}\label{f6.29}
\begin{aligned}
\bar{e}_{\pm, n}^{(1)}&= {\cal E}(V^{\pm, n+1}, \Phi^{\pm,
n+1})-{\cal E}(V^{\pm, n}, \Phi^{\pm, n}) -{\cal E}'_{(V^{\pm, n},
\Phi^{\pm,n})}(\delta {V}^{\pm, n}, \delta \Phi^{\pm,n})\cr
&=\partial_{x_2}(\delta \Phi^{\pm,n})\delta V_3^{\pm,n}+
\partial_{x_3}(\delta \Phi^{\pm,n})\delta V_4^{\pm,n},
\end{aligned}
\end{equation}
the second error resulting from the substitution:
 %\stepcounter{nonl}
\begin{equation}\label{f6.30}
\begin{aligned}
\bar{e}_{\pm,n}^{(2)}=&{\cal E}'_{(V^{\pm, n}, \Phi^{\pm,
n})}(\delta {V^{\pm, n}},
 \delta {\Phi^{\pm,n}})
-{\cal E}'_{(S_{\theta_n}V^{\pm, n}, S_{\theta_n}\Phi^{\pm,n})}
(\delta {V^{\pm, n}},
 \delta {\Phi^{\pm,n}})\cr
=&\partial_{x_2}((I-S_{\theta_n})\Phi^{\pm,n})\delta V_3^{\pm,n}
+\partial_{x_3}((I-S_{\theta_n})\Phi^{\pm,n})\delta V_4^{\pm,n}\cr
&+(I-S_{\theta_n})V_3^{\pm,n}\partial_{x_2}(\delta\Phi^{\pm,n})
+(I-S_{\theta_n})V_4^{\pm,n}\partial_{x_3}(\delta\Phi^{\pm,n}),
\end{aligned}
\end{equation}
the third error resulting from the second substitution:
%\stepcounter{nonl}
\begin{equation}\label{f6.31}
\bar{e}_{\pm,n}^{(3)}= {\cal E}'_{(S_{\theta_n}V^{\pm, n},
S_{\theta_n}\Phi^{\pm,n})}(\delta {V^{\pm, n}},
 \delta {\Phi^{\pm,n}})
-{\cal E}'_{(V^{\pm, n+\frac 12}, S_{\theta_n}\Phi^{\pm,n})}(\delta
{V^{\pm, n}},
 \delta {\Phi^{\pm,n}})
\end{equation}
which vanishes due to (\ref{f6.24}) and that all the coefficients of
the linearized operator ${\cal E}_{(V^\pm, \Phi^\pm)}'$ are
independent of $(V_2^\pm, V_5^\pm)$, and the remaining error:
%\stepcounter{nonl}
\begin{equation}\label{f6.32}
\bar{e}_{\pm,n}^{(4)}= {\cal E}'_{(V^{\pm, n+\frac 12},
S_{\theta_n}\Phi^{\pm,n})}(\delta {V^{\pm, n}},
 \delta {\Phi^{\pm,n}})
-{\cal E}'_{(V^{\pm, n+\frac 12}, S_{\theta_n}\Phi^{\pm,n})}(\delta
{\dot{V}^{\pm, n}},
 \delta {\Phi^{\pm,n}}).
\end{equation}

Thus, we have
%\stepcounter{nonl}
\begin{equation}\label{f6.33}
{\cal E}(V^{\pm, n+1}, \Phi^{\pm,n+1})=\sum_{j=0}^n ({\cal
E}'_{(V^{\pm, j+\frac 12}, S_{\theta_j}\Phi^{\pm,j})}(\delta
{V^{\pm, j}},
 \delta {\Phi^{\pm,j}})+\bar{e}_{\pm, j})
\end{equation}
which leads to define $\delta\Phi^{\pm, n}$ that is governed by the
following problem:
%\stepcounter{nonl}
\begin{equation}\label{f6.34}
\begin{cases}
{\cal E}'_{(V^{\pm, n+\frac 12}, S_{\theta_n}\Phi^{\pm,n})}(\delta
\dot{V}^{\pm, n},
 \delta {\Phi^{\pm,n}})=h_n^\pm \qquad {\rm on}~ \Omega_T,\cr
\delta {\Phi^{\pm,n}}|_{t\le 0}=0,
\end{cases}
\end{equation}
where  $h_n^\pm$ are defined by
%\stepcounter{nonl}
\begin{equation}\label{f6.35}
\sum_{k=0}^nh_k^\pm+S_{\theta_n}\big(\sum_{k=0}^{n-1}\bar{e}_{\pm,
k}\big)=0
\end{equation}
by induction on $n$ starting with $h_0^\pm=0$.

By comparing (\ref{f6.34}) with (\ref{f6.21}), it is easy to verify
%\stepcounter{nonl}
\begin{equation}\label{f6.36}
\delta {\Phi^{\pm,n}}|_{x_1=0}=\delta \phi^{n}.
\end{equation}

\bigskip
{}From \eqref{f5.11}--\eqref{f5.12} and \eqref{f6.34}, we know that
the compatibility conditions hold for the initial-boundary value
problem \eqref{f6.21} for all $n\ge 0$.

\medskip
The steps for determining $(\delta{V}^{\pm, n}, \delta\Phi^{\pm,
n})$ are to solve first $\delta\dot{V}^{\pm, n}$ from (\ref{f6.21})
and then $\delta\Phi^{\pm, n}$ from (\ref{f6.34}). Finally,
$\delta{V}^{\pm, n}$ can be obtained from (\ref{f6.8}).

\section{Nonlinear Stability III: Convergence of the Iteration Scheme\\[1mm]
and Existence of the Current-Vortex Sheet}\label{s:7}

Fix any $s_0\ge 9$, $\alpha\ge s_0+5$, and $s_1\in
[\alpha+5,2\alpha-s_0]$. Let the zero-th order approximate solutions
for the initial data $(U_0^\pm, \psi_0)$ constructed in \S5 satisfy
%\stepcounter{nonl}
\begin{equation}\label{f7.1}
\|\dot{U}_a^\pm\|_{s_1+3, T}+\|\dot{\Psi}_a^\pm\|_{s_1+3, T}
+\|f_a^\pm\|_{s_1-4, T}\le \delta, \qquad \|f_a^\pm\|_{\alpha+1,
T}/\delta \quad {\rm is~small}
\end{equation}
for some small constant $\delta>0$.

Before we prove the convergence of the iteration scheme
(\ref{f6.21}) and (\ref{f6.34}), we first introduce the following
lemmas.

\vspace{.1in} {\bf Lemma 7.1.} {\it For any $k\in [0, n]$,}
%\stepcounter{nonl}
\begin{equation}\label{f7.13}
\begin{cases} \|(V^{\pm,k},
\Phi^{\pm,k})\|_{s,T}+\|\phi^{k}\|_{H^{s-1}(b\Omega_T)} \le
C\delta\, \theta_k^{(s-\alpha)_+} \quad \mbox{for}\,\, s\in [s_0,
s_1], ~s\neq \alpha,\cr \|(V^{\pm,k},
\Phi^{\pm,k})\|_{\alpha,T}+\|\phi^{k}\|_{H^{\alpha-1}(b\Omega_T)}
\le C\delta \log\theta_k,\cr \|(S_{\theta_k}V^{\pm,k},
S_{\theta_k}\Phi^{\pm,k})\|_{s,T}
+\|S_{\theta_k}\phi^{k}\|_{H^{s-1}(b\Omega_T)} \le C\delta
\,\theta_k^{(s-\alpha)_+} \quad\mbox{for}\,\, s\ge s_0, s\neq
\alpha,\cr \|(S_{\theta_k}V^{\pm,k},
S_{\theta_k}\Phi^{\pm,k})\|_{\alpha,T}
+\|S_{\theta_k}\phi^{k}\|_{H^{\alpha-1}(b\Omega_T)} \le C\delta
\log\theta_k,\cr \|((I-S_{\theta_k})V^{\pm,k},
(I-S_{\theta_k})\Phi^{\pm,k})\|_{s,T}
+\|(I-S_{\theta_k})\phi^{k}\|_{H^{s-1}(b\Omega_T)} \le C\delta
\,\theta_k^{s-\alpha} \quad\mbox{for}\,\, s\in [s_0, s_1].
\end{cases}
\end{equation}
These results can be easily obtained by using the triangle
inequality, the classical comparison between series and integrals,
and the properties in \eqref{f6.2} of the smoothing operators
$S_\theta$.

\vspace{.1in} {\bf Lemma 7.2.} {\it If $4\le s_0\le \alpha$ and
$\alpha\ge 5$, then, for the modified state $V^{\pm, n+\frac 12}$,
we have
%\stepcounter{nonl}
\begin{equation}\label{f7.14}
\|V^{\pm, n+\frac 12}-S_{\theta_n}V^{\pm, n}\|_{s, T}\le C\delta
\theta_n^{s+1-\alpha}\qquad\mbox{for any}\,\, s\in [4, s_1+2].
\end{equation} }

{\bf Lemma 7.3.} {\it Let $s_0\ge 5$ and $\alpha\ge s_0+ 3$. For all
$k\le n-1$, we have
%\stepcounter{nonl}
\begin{equation}\label{f7.15}
\begin{cases} \|e_{\pm,k}\|_{s,T}\le
C\delta^2\theta_k^{L^1(s)}\Delta_k \qquad\mbox{for}\,\, s\in[s_0,
s_1-4],\cr \|\bar{e}_{\pm,k}\|_ {s,T}\le
C\delta^2\theta_k^{{L^2(s)}}\Delta_k \qquad\mbox{for}\,\, s\in [s_0,
s_1-3],\cr \|\tilde{e}_{k}\|_{H^s(b\Omega_T)}\le
C\delta^2\theta_k^{{L^3(s)}}\Delta_k \qquad\mbox{for}\,\, s\in
[s_0-1, s_1-4],
\end{cases}
\end{equation}
where
%\stepcounter{nonl}
\begin{eqnarray*}
&&L^1(s)=
\begin{cases}
\max((s+2-\alpha)_++s_0-\alpha-1, s+s_0+4-2\alpha)\qquad {\rm
if}\quad \alpha\neq s+2, s+4,\cr \max(s_0-\alpha,
2(s_0-\alpha)+4)\qquad {\rm if}\quad \alpha=s+4,\cr
s_0+2-\alpha\qquad {\rm if}\quad \alpha=s+2,\end{cases}\\
&& L^2(s)=
\begin{cases}
\max((s+2-\alpha)_++2(s_0-\alpha), s+s_0+2-2\alpha)\qquad {\rm
if}\quad \alpha\neq s+2, s+3,\cr s_0-\alpha-1\qquad {\rm if}\quad
\alpha=s+3,\cr s_0-\alpha\qquad {\rm if}\quad \alpha=s+2,
\end{cases}
\end{eqnarray*}
and
$$
L^3(s)=
\begin{cases}
\max((s+3-\alpha)_++2(s_0-\alpha), s+s_0+3-2\alpha)\qquad {\rm
if}\quad \alpha\neq s+3, s+4,\cr s_0-\alpha-1\qquad {\rm if}\quad
\alpha=s+4,\cr s_0-\alpha\qquad {\rm if}\quad \alpha=s+3.
\end{cases}
$$
}

Denote the accumulated errors by
%\stepcounter{nonl}
\begin{equation}\label{f7.17}
E_{\pm, n}=\sum_{k=0}^{n-1}e_{\pm,k}, \quad
\tilde{E}_{n}=\sum_{k=0}^{n-1}\tilde{e}_{k}, \quad \bar{E}_{\pm,
n}=\sum_{k=0}^{n-1}\bar{e}_{\pm,k}.
\end{equation}
Then, as a corollary of Lemma 7.3, we have

\medskip
{\bf Lemma 7.4.} {\it Let $s_0\ge 5$, $\alpha\ge s_0+3$, and $s_1\le
2\alpha-s_0-1$. Then
%\stepcounter{nonl}
\begin{equation}\label{f7.18}
\begin{cases} \|E_{\pm,n}\|_{s,T}\le C\delta^2\theta_n \qquad\mbox{for}\,\,
s\in [s_0, s_1-4],\cr \|\bar{E}_{\pm,n}\|_ {s,T}\le \delta^2
\qquad\mbox{for}\,\, s\in [s_0, s_1-3],\cr
\|\tilde{E}_{n}\|_{H^s(b\Omega_T)}\le C\delta^2\theta_n
\qquad\mbox{for}\,\, s\in [s_0-1, s_1-4].
\end{cases}
\end{equation}}

\vspace{.1in} {\bf Lemma 7.5.} {\it For any $s_0\ge 5$, $\alpha\ge
s_0+5$, and $s_1\in [\alpha+5, 2\alpha-s_0-1]$, we have
%\stepcounter{nonl}
\begin{equation}\label{f7.19}
\begin{cases} \|f_n^\pm\|_{s,T}\le C\Delta_n
\theta_n^{s-\alpha-1}(\|f_a^\pm\|_{\alpha,T}+\delta^2),
%\qquad%\forall \tilde{s}\ge 0
\cr \|g_n\|_{H^{s+1}(b\Omega_T)}\le C\delta^2\Delta_n
\theta_n^{s-\alpha-1},\cr \|h^\pm_n\|_{s,T}\le C\delta^2\Delta_n
\theta_n^{s-\alpha-1}
\end{cases}
\end{equation}
for all $s\ge s_0$. }

\medskip
The proofs of Lemmas 7.2--7.3 and 7.5 will be given in Section 8.
With these lemmas, we can now prove the following key result for the
convergence of the iteration scheme.

\vspace{.1in} {\bf Proposition 7.1.} {\it For the solution sequence
$(\delta V^{\pm k}, \delta \Phi^{\pm, k}, \delta \phi^k)$ given by
{\rm (\ref{f6.21})} and {\rm (\ref{f6.34})}, we have
\begin{equation}\label{f7.2}
\begin{cases}
\|(\delta V^{\pm,k}, \delta \Phi^{\pm,k})\|_{s,T}+\|\delta
\phi^{k}\|_{H^{s-1}(b\Omega_T)} \le\delta
\theta_k^{s-\alpha-1}\Delta_k \qquad \mbox{for}\,\,
s\in[s_0,s_1],\cr \|{\cal L}(V^{\pm,k+1},
\Phi^{\pm,k+1})V^{\pm,k+1}-f_a^\pm\|_{s,T} \le 2\delta
\theta_{k+1}^{s-\alpha-1} \qquad\mbox{for}\,\, s\in [s_0, s_1-4],\cr
\|{\cal B}(V^{+,k+1}, V^{-,k+1}, \phi^{k+1})\|_{H^{s-1}(b\Omega_T)}
\le \delta \theta_{k+1}^{s-\alpha-1} \qquad\mbox{for}\,\, s\in [s_0,
s_1-2]
\end{cases}
\end{equation}
for any $k\ge 0$, where $\Delta_k=\theta_{k+1}-\theta_k$.}

\medskip
{\it Proof:} Estimate (\ref{f7.2}) is proved by induction on $k\ge
0$.

{\it Step 1. Verification of {\rm (\ref{f7.2})} for $k=0$.} We first
notice from (\ref{f5.9})--(\ref{f5.10}) that $(U_a^\pm, \Psi_a^\pm)$
satisfies the Rankine-Hugoniot conditions, and $V^{\pm,
0}=\Phi^{\pm, 0}=0$ imply that $V^{\pm, \frac 12}=0$. Thus,
$\delta\dot{V}^{\pm, 0}$ satisfies the following problem:
%\stepcounter{nonl}
\begin{equation}\label{f7.3}
\begin{cases}
L'_{e,(U_a^\pm, \Psi_a^\pm)}\delta \dot{V}^{\pm, 0}
=S_{\theta_0}f_a^\pm \qquad {\rm in}~\Omega_T,\cr B'_{(U_a^\pm,
\psi_a)}(\delta \dot{V}^{+, 0}, \delta \dot{V}^{-, 0},
\delta\phi^0)=0 \qquad {\rm on}~b\Omega_T,\cr
\delta \dot{V}^{\pm,
0}|_{t\le 0}=0.
\end{cases}
\end{equation}

Applying Theorem 4.2 to problem (\ref{f7.3}), we have
%\stepcounter{nonl}
\begin{equation}\label{f7.4}
\|\delta \dot{V}^{\pm, 0}\|_{s,T}\le
C\theta_0^{(s-\alpha-1)_+}\|f_a^\pm\|_{\alpha+1, T} \qquad\mbox{for
all}\,\, s\in [s_0, s_1].
\end{equation}

Similarly, from (\ref{f6.34}), we obtain
%\stepcounter{nonl}
\begin{equation}\label{f7.5}
\|\delta \Phi^{\pm, 0}\|_{s,T}\le
C\left(\|\tilde{h}_0^\pm\|_{s,T}+\|\dot{U}_a^\pm\|_{s,T}\|\tilde{h}_0^\pm\|_{s_0,T}\right),
\end{equation}
where $\tilde{h}_0^\pm=\delta \dot{V_2}^{\pm,
0}-\partial_{x_2}\Psi_a^\pm \delta \dot{V_3}^{\pm, 0}
-\partial_{x_3}\Psi_a^\pm \delta \dot{V_4}^{\pm, 0}$ satisfies
\begin{equation}\label{f7.6}
\|\tilde{h}_0^\pm\|_{s,T}\le
C\|\delta\dot{V}^{\pm,0}\|_{s,T}\le
C_0\theta_0^{(s-\alpha-1)_+}\|f_a^\pm\|_{\alpha+1, T}.
\end{equation}
Thus, we have
\begin{equation}\label{f7.7}
\|\delta \Phi^{\pm, 0}\|_{s,T}\le
C_1\theta_0^{(s-\alpha-1)_+}\|f_a^\pm\|_{\alpha+1, T}.
\end{equation}

Since
$$
\delta\dot{V}^{\pm,0}=\delta{V}^{\pm,0}
-\delta\Phi^{\pm,0}\frac{\partial_{x_1}\dot{U}_a^\pm}{\partial_{x_1}\Psi_a^\pm},
$$
we use (\ref{f7.4}) and (\ref{f7.7}) to find
\begin{equation}\label{f7.8}
\|\delta V^{\pm, 0}\|_{s,T}\le
C_2\theta_0^{(s-\alpha-1)_+}\|f_a^\pm\|_{\alpha+1, T}.
\end{equation}

{}From (\ref{f7.7})--(\ref{f7.8}), we deduce
\begin{equation}\label{f7.9}
\|(\delta V^{\pm, 0}, \delta\Phi^{\pm,0})\|_{s,T}\le
\delta\theta_0^{s-\alpha-1}\Delta_0\qquad \mbox{for any}\,\, s\in
[s_0, s_1],
\end{equation}
provided that $\|f_a^\pm\|_{\alpha+1, T}/\delta$ is small.

Obviously, one has
\begin{equation}\label{f7.10}
\|{\cal L}(V^{\pm, 1}, \Phi^{\pm,1})V^{\pm, 1}-f_a^\pm\|_{s,T}\le
\|(S_{\theta_0}-I)f_a^\pm\|_{s,T}+\|e_{\pm, 0}\|_{s,T},
\end{equation}
where
$$
\begin{array}{ll}e_{\pm, 0}
= & L(U_a^\pm+V^{\pm, 1}, \Psi_a^\pm+\Phi^{\pm, 1})(U_a^\pm+V^{\pm, 1})
-L(U_a^\pm, \Psi_a^\pm)U_a^\pm\\
&
-L'_{(U_a^\pm, \Psi_a^\pm)}(\delta V^{\pm, 0}, \delta\Phi^{\pm, 0})-
\frac{\delta\Phi^{\pm, 0}}{\partial_{x_1}\Psi_a^\pm}\partial_{x_1}f_a^\pm\end{array}
$$
satisfies
\begin{equation}\label{f7.11}
\|e_{\pm, 0}\|_{s,T}\le
C\theta_0^{(s-\alpha-1)_+}\|f_a^\pm\|_{\alpha+1, T}
\qquad\mbox{for}\,\, s\in [s_0, s_1-2].
\end{equation}

Using the property of $S_{\theta_0}$, we have
$$
\|(S_{\theta_0}-I)f_a^\pm\|_{s,T} \le
C\theta_0^{s-\alpha-1}\|f_a^\pm\|_{\alpha+1, T} \qquad
\mbox{for}\,\, s\le \alpha+1,
$$
and
$$
\|(S_{\theta_0}-I)f_a^\pm\|_{s,T}\le
C\theta_0^{s-s_1}\|f_a^\pm\|_{s_1, T} \qquad \mbox{for}\,\, s\in
[\alpha+1, s_1].
$$

Thus, from (\ref{f7.10}), we obtain
\begin{equation}\label{f7.12}
\|{\cal L}(V^{\pm, 1}, \Phi^{\pm,1})V^{\pm, 1}-f_a^\pm\|_{s,T}\le
2\delta\theta_1^{s-\alpha-1}\qquad\mbox{for}\,\, s\in [s_0,s_1-2],
\end{equation}
when $\|f_a^\pm\|_{\alpha+1, T}/\delta$ is small and $\theta_0>0$ is
properly large.

Similarly, we can obtain
$$
\|{\cal B}(V^{+,1}, V^{-,1}, \phi^{1})\|_{H^{s-1}(b\Omega_T)} \le
\delta \theta_{1}^{s-\alpha-1} \qquad\mbox{for}\,\, s\in
[s_0,s_1-2].
$$

\medskip
{\it Step 2. Suppose that {\rm (\ref{f7.2})} holds for all $k\le
n-1$, we now verify {\rm (\ref{f7.2})} for $k=n$.} First, we note
that, when $\alpha\ge s_0+3$,
\begin{equation}\label{f7.20}
\begin{aligned}
&\|(\dot{U}_a^\pm+V^{\pm, n+\frac 12},
\dot{\Psi}_a^\pm+S_{\theta_n}\Phi^{\pm, n})\|_{s_0+2, T}\\
&\le \|(\dot{U}_a^\pm, \dot{\Psi}_a^\pm)\|_{s_0+2, T} +\|V^{\pm,
n+\frac 12}-S_{\theta_n}V^{\pm, n}\|_{s_0+2, T}+
\|(S_{\theta_n}V^{\pm, n},S_{\theta_n}\Phi^{\pm, n})\|_{s_0+2, T}\cr
&\le C\delta(1+\theta_n^{s_0+3-\alpha}+\theta_n^{(s_0+2-\alpha)_+})
\le C\delta,
\end{aligned}
\end{equation}
by using assumption (\ref{f7.1}) and Lemmas 7.1--7.2. Applying
Theorem 4.2 to problem (\ref{f6.21}), we have
%\stepcounter{nonl}
\begin{equation}\label{f7.21}
\begin{aligned}
&\|\delta\dot{V}^{\pm, n}\|_{s,T}\le
C\big(\|f_n^\pm\|_{s,T}+\|g_n\|_{H^{s+1}(b\Omega_T)}
+\|\dot{coef}\|_{s,T}(\|f_n^\pm\|_{s_0,T}+\|g_n\|_{H^{s_0+1}(b\Omega_T)})\big)\cr
&\le C\big(\|f_n^\pm\|_{s,T}+\|g_n\|_{H^{s+1}(b\Omega_T)}
+\|(\dot{U}_a^\pm+V^{\pm, n+\frac 12},
\dot{\Psi}_a^\pm+S_{\theta_n}\Phi^{\pm,
n})\|_{s+2,T}(\|f_n^\pm\|_{s_0,T}+\|g_n\|_{H^{s_0+1}(b\Omega_T)})\big).
\end{aligned}
\end{equation}

On the other hand, similar to \eqref{f7.20}, from assumption
(\ref{f7.1}) and Lemmas 7.1--7.2, we have
%\stepcounter{nonl}
\begin{equation}\label{f7.22}
\|(\dot{U}_a^\pm+V^{\pm, n+\frac 12},
\dot{\Psi}_a^\pm+S_{\theta_n}\Phi^{\pm, n})\|_{s+2,T}\le
C\delta\big(1+\theta_n^{s+3-\alpha}+\theta_n^{(s+2-\alpha)_+}\big),
\end{equation}
which implies
%\stepcounter{nonl}
\begin{equation}\label{f7.23}
\|(\dot{U}_a^\pm+V^{\pm, n+\frac 12},
\dot{\Psi}_a^\pm+S_{\theta_n}\Phi^{\pm, n})\|_{s+2,T}
(\|f_n^\pm\|_{s_0,T}+\|g_n\|_{H^{s_0+1}(b\Omega_T)})\le
C\theta_n^{s-\alpha-1}\Delta_n(\|f_a^\pm\|_{\alpha, T}+\delta^2)
\end{equation}
for all $s\in [s_0, s_1]$ by using $\alpha\ge s_0+3$ and
(\ref{f7.19}).

Thus, from (\ref{f7.21}), we conclude
%\stepcounter{nonl}
\begin{equation}\label{f7.24}
\|\delta \dot{V}^{\pm, n}\|_{s,T} \le
C\theta_n^{s-\alpha-1}\Delta_n(\|f_a^\pm\|_{\alpha, T}+\delta^2)
\qquad\mbox{for}\,\, s\in [s_0, s_1].
\end{equation}

For problem (\ref{f6.34}), we can easily obtain the following estimate:
\begin{equation}\label{f7.25}
\|\delta \Phi^{\pm, n}\|_{s,T} \le
C(\|\tilde{h}_n^\pm\|_{s,T}+\|\dot{U}_a^\pm+S_{\theta_n}V^{\pm,
n}\|_{s,T} \|\tilde{h}_n^\pm\|_{s_0,T}) \qquad\mbox{for any}\,\,
s\ge s_0,
\end{equation}
where $C>0$ depends only on
$\|\dot{U}_a^\pm+S_{\theta_n}V^{\pm, n}\|_{s_0,T}$, and
\begin{equation}\label{f7.26}
\tilde{h}_n^\pm=h_n^\pm+\delta \dot{V}_2^{\pm, n}
-\partial_{x_2}(\dot{\Psi}_a^\pm+S_{\theta_n}\Phi^{\pm, n})\delta \dot{V}_3^{\pm, n}
-\partial_{x_3}(\dot{\Psi}_a^\pm+S_{\theta_n}\Phi^{\pm, n})\delta \dot{V}_4^{\pm, n}.
\end{equation}

Obviously, we have
$$
\|\tilde{h}_n^\pm\|_{s,T}\le \|{h}_n^\pm\|_{s,T}
+(1+\|\dot{\Psi}_a^\pm+S_{\theta_n}\Phi^{\pm,
n}\|_{s_0,T})\|\delta\dot{V}^{\pm, n}\|_{s,T}
+\|\dot{\Psi}_a^\pm+S_{\theta_n}\Phi^{\pm,
n}\|_{s+1,T}\|\delta\dot{V}^{\pm, n}\|_{s_0,T},
$$
which implies
\begin{equation}\label{f7.27}
\|\tilde{h}_n^\pm\|_{s,T}\le C\Delta_n\theta_n^{s-\alpha-1}(\|f_a^\pm\|_{\alpha, T}+\delta^2)
\end{equation}
by using (\ref{f7.24}) and Lemmas 7.1 and 7.5.

\medskip
Substituting (\ref{f7.27}) into (\ref{f7.25}), we find
\begin{equation}\label{f7.28}
\|\delta \Phi^{\pm, n}\|_{s,T} \le
C\Delta_n\theta_n^{s-\alpha-1}(\|f_a^\pm\|_{\alpha, T}+\delta^2)
\qquad\mbox{for all}\,\, s\in [s_0, s_1].
\end{equation}

\medskip
Together (\ref{f7.24}) with (\ref{f7.28}), we obtain (\ref{f7.2})
for $(\delta V^{\pm, n},\delta \Phi^{\pm,n})$ by using
$$
\delta V^{\pm, n}
=\delta\dot{V}^{\pm,n}+\frac{\partial_{x_1}(U_a^\pm+V^{\pm, n+\frac 12})}
{\partial_{x_1}(\Psi_a^\pm+ S_{\theta_n}\Phi^{\pm, n})}\delta\Phi^{\pm, n}
$$
and choosing $\|f_a^\pm\|_{\alpha, T}/\delta$ and $\delta>0$ to be
small.

\medskip
Finally, we
verify the other inequalities in (\ref{f7.2}). First, from
(\ref{f6.21})--(\ref{f6.22}), we have
%\stepcounter{nonl}
\begin{equation}\label{f7.29}
{\cal L}(V^{\pm,n}, \Phi^{\pm,n})V^{\pm,n}-f_a^\pm
=(S_{\theta_{n-1}}-I)f_a^\pm+(I-S_{\theta_{n-1}})E_{\pm, n-1}+
e_{\pm, n-1}.
\end{equation}

{}From Lemma 7.4, we have
%\stepcounter{nonl}
\begin{equation}\label{f7.30}
\|(I-S_{\theta_{n-1}})E_{\pm, n-1}\|_{s,T}
\le C\theta_n^{s-\tilde{s}}\|E_{\pm, n-1}\|_{\tilde{s},T}
\le C\theta_n^{s-\tilde{s}+1}\delta^2\le C\theta_n^{s-\alpha-1}\delta^2
\end{equation}
by choosing $\tilde{s}=\alpha+2$.

{}From Lemma 7.3, we obtain
%\stepcounter{nonl}
\begin{equation}\label{f7.31}
\|e_{\pm, n-1}\|_{s,T} \le C\theta_n^{L^1(s)}\delta^2\le
C\theta_n^{s-\alpha-1}\delta^2 \qquad\mbox{for all}\,\, s\in [s_0,
s_1-4],
\end{equation}
by using
$\alpha\ge s_0+5$.

When $s\le \alpha+1$, we have
%\stepcounter{nonl}
\begin{equation}\label{f7.32}
\|(S_{\theta_{n-1}}-I)f_a^\pm\|_{s,T} \le
C\theta_n^{s-\alpha-1}\|f_a^\pm\|_{\alpha+1, T};
\end{equation}
while $s\in (\alpha+1, s_1-4]$, we have
%\stepcounter{nonl}
\begin{equation}\label{f7.33}
\|(S_{\theta_{n-1}}-I)f_a^\pm\|_{s,T}\le
\|S_{\theta_{n-1}}f_a^\pm\|_{s,T}+\| f_a^\pm\|_{s,T}\le
C\theta_n^{s-\alpha-1}\|f_a^\pm\|_{\alpha+1, T}+\delta
\end{equation}
by using (\ref{f6.2}) and (\ref{f7.1}).

Substitution \eqref{f7.30}--\eqref{f7.33} into (\ref{f7.29}) yields
%\stepcounter{nonl}
\begin{equation*}\label{f7.34}
\|{\cal L}(V^{\pm,n}, \Phi^{\pm,n})V^{\pm,n}-f_a^\pm\|_{s,T} \le
2\delta \theta_n^{s-\alpha-1} \qquad\mbox{for}\,\, s\in [s_0,
s_1-4],
\end{equation*}
provided that
%\stepcounter{nonl}
%\begin{equation}\label{f7.35}
$\|f_a^\pm\|_{s_1-4,T} \le \delta$ and
$\|f_a^\pm\|_{\alpha+1,T}/\delta\quad {\rm is ~small}.$
%\end{equation}

Similarly, we have
%\stepcounter{nonl}
\begin{equation}\label{f7.36}
\|{\cal B}(V^{+,n},V^{-,n}, \phi^{n})\|_{H^{s-1}(b\Omega_T)}
\le \delta \theta_n^{s-\alpha-1}.
\end{equation}
Thus, we obtain (\ref{f7.2}) for $k=n$.

\vspace{.2in} {\sc Convergence of the Iteration Scheme:} {}From
Proposition 7.1, we have
%\stepcounter{nonl}
\begin{equation}\label{f7.37}
\sum_{n\ge 0}\big(\|(\delta V^{\pm, n}, \delta\Phi^{\pm,
n})\|_{\alpha,T}+\|\delta\phi^n\|_{H^{\alpha-1}(b\Omega_T)}\big)<\infty,
\end{equation}
which implies that there exists $(V^\pm, \Phi^\pm)\in
B^{\alpha}(\Omega_T)$ with $\phi\in H^{\alpha-1}(b\Omega_T)$ such
that
%\stepcounter{nonl}
\begin{equation}\label{f7.38}
\begin{cases}(V^{\pm, n}, \Phi^{\pm, n})\longrightarrow (V^{\pm},
\Phi^{\pm})\qquad {\rm in}\,\, B^{\alpha}(\Omega_T)\times
B^{\alpha}(\Omega_T),\cr \phi^n\longrightarrow \phi\quad {\rm
in}\quad H^{\alpha-1}(b\Omega_T).
\end{cases}
\end{equation}

Thus, we conclude
%the Main Theorem in Section 2. \iffalse

%stop

\vspace{.1in} {\bf Theorem 7.1.} {\it Let $\alpha\ge 14$ and $s_1\in
[\alpha+5, 2\alpha-9]$. Let $\psi_0\in H^{2s_1+3}(\R^2)$ and
$U_0^\pm-\bar{U}^\pm\in B^{2(s_1+2)}(\R^3_+)$ satisfy the
compatibility conditions of problem \eqref{f2.15}--\eqref{f2.18} up
to order $s_1+2$, and let conditions \eqref{f2.3}--\eqref{f2.3a} and
\eqref{f7.1} be satisfied. Then there exists a solution $(V^\pm,
\Phi^\pm)\in B^\alpha(\Omega_T)$ with $\phi\in
H^{\alpha-1}(b\Omega_T)$ to problem \eqref{f5.14}.}

Then Theorem 2.1 (Main Theorem) in Section 2 directly follows from
Theorem 7.1.

\section{Error Estimates: Proofs of Lemmas 7.2--7.3 and 7.5}

In this section, we study the error estimates for the iteration
scheme (\ref{f6.21})--(\ref{f6.22}) and (\ref{f6.34})--(\ref{f6.35})
to provide the proofs for Lemmas 7.2--7.3 and 7.5 under the
assumption that (\ref{f7.2}) holds for all $0\le k\le n-1$. We start
with the proof of Lemma 7.2.

\vspace{.1in} {\sc Proof of Lemma {\rm 7.2.}} Denote by
%\stepcounter{nonl}
\begin{equation}\label{f8.1}
{\cal E}_1^{\pm, n}:={\cal E}(V^{\pm, n}, \Phi^{\pm, n}),
\end{equation}
and
%\stepcounter{nonl}
\begin{equation}\label{f8.2}
%\begin{aligned}
{\cal E}_2^{\pm, n}:= V_5^{\pm, n}-
\partial_{x_2}(\Psi^{\pm}_a+\Phi^{\pm, n})V_6^{\pm, n}
-\partial_{x_3}(\Psi^{\pm}_a+\Phi^{\pm, n})V_7^{\pm,
n}-U^\pm_{a,6}\partial_{x_2}\Phi^{\pm,
n}-U^\pm_{a,7}\partial_{x_3}\Phi^{\pm, n}
%\end{aligned}
\end{equation}
as the extension of $({\cal B}(V^{\pm, n}, \phi^{n}))_2^\pm$ in
$\Omega_T$.

By the definition in (\ref{f6.25}) for $V^{\pm, n+\frac 12}$, we
have
%\stepcounter{nonl}
\begin{equation}\label{f8.3}\begin{cases}
V_2^{\pm, n+\frac 12}-S_{\theta_n}V_2^{\pm, n}\\
\hspace{.3in}=S_{\theta_n}{\cal E}_1^{\pm, n}+[\partial_t,
S_{\theta_n}]\Phi^{\pm,
n}+
%\hspace{.3in}
\big(\partial_{x_2}(\Psi^{\pm}_a+S_{\theta_n}\Phi^{\pm,
n})S_{\theta_n}V_3^{\pm, n}-
S_{\theta_n}(\partial_{x_2}(\Psi^{\pm}_a+\Phi^{\pm, n})V_3^{\pm,
n})\big)\cr \hspace{.5in}+
\big(\partial_{x_3}(\Psi^{\pm}_a+S_{\theta_n}\Phi^{\pm,
n})S_{\theta_n}V_4^{\pm, n}-
S_{\theta_n}(\partial_{x_3}(\Psi^{\pm}_a+\Phi^{\pm, n})V_4^{\pm,
n})\big)\cr \hspace{.5in}+\big(\partial_{x_2}(S_{\theta_n}\Phi^{\pm,
n})U_{a,3}^{\pm}- S_{\theta_n}(\partial_{x_2}\Phi^{\pm,
n}U_{a,3}^{\pm})\big)+ \big(\partial_{x_3}(S_{\theta_n}\Phi^{\pm,
n})U_{a,4}^{\pm}-
S_{\theta_n}(\partial_{x_3}\Phi^{\pm,n}U_{a,4}^{\pm})\big),\cr
V_5^{\pm,
n+\frac 12}-S_{\theta_n}V_5^{\pm, n}\\
\hspace{.3in}=-S_{\theta_n}{\cal E}_2^{\pm, n}
%\hspace{.3in}
+\big(\partial_{x_2}(\Psi^{\pm}_a+S_{\theta_n}\Phi^{\pm,
n})S_{\theta_n}V_6^{\pm, n}-
S_{\theta_n}(\partial_{x_2}(\Psi^{\pm}_a+\Phi^{\pm, n})V_6^{\pm,
n})\big)\cr \hspace{.5in}+
\big(\partial_{x_3}(\Psi^{\pm}_a+S_{\theta_n}\Phi^{\pm,
n})S_{\theta_n}V_7^{\pm, n}-
S_{\theta_n}(\partial_{x_3}(\Psi^{\pm}_a+\Phi^{\pm, n})V_7^{\pm,
n})\big)\cr \hspace{.5in}+\big(\partial_{x_2}(S_{\theta_n}\Phi^{\pm,
n})U_{a,6}^{\pm}- S_{\theta_n}(\partial_{x_2}\Phi^{\pm,
n}U_{a,6}^{\pm}))+ (\partial_{x_3}(S_{\theta_n}\Phi^{\pm,
n})U_{a,7}^{\pm}-
S_{\theta_n}(\partial_{x_3}\Phi^{\pm,n}U_{a,7}^{\pm})\big).
\end{cases}
\end{equation}

On the other hand, we have
%\stepcounter{nonl}
\begin{equation}\label{f8.4}\begin{cases}
{\cal E}_1^{\pm, n}={\cal E}_1^{\pm, n-1}+
%(V^{\pm, n-1}, \Phi^{\pm,n-1})+
\partial_t(\delta\Phi^{\pm, n-1})
+\partial_{x_2}(\Psi^{\pm}_a+\Phi^{\pm, n}) \delta V_3^{\pm, n-1}\cr
\hspace{.5in} +\partial_{x_3}(\Psi^{\pm}_a+\Phi^{\pm, n})\delta
V_4^{\pm, n}-\delta V_2^{\pm, n-1}+
\partial_{x_2}(\delta\Phi^{\pm,n-1})(U_{a,3}^\pm+V_3^{\pm, n-1})
\cr\hspace{.5in}
+\partial_{x_3}(\delta\Phi^{\pm,n-1})(U_{a,4}^\pm+V_4^{\pm,
n-1}),\cr {\cal E}_2^{\pm, n}= {\cal E}_2^{\pm, n-1}
% {\cal E}_T({\cal B}(V^{\pm,n-1},\phi^{n-1})))_2^\pm
-\partial_{x_2}(\Psi^{\pm}_a+\Phi^{\pm, n}) \delta V_6^{\pm, n-1}\cr
\hspace{.5in} -\partial_{x_3}(\Psi^{\pm}_a+\Phi^{\pm, n})\delta
V_7^{\pm, n}+\delta V_5^{\pm, n-1}-
\partial_{x_2}(\delta\Phi^{\pm,n-1})(U_{a,6}^\pm+V_6^{\pm, n-1})
\cr\hspace{.5in}
+\partial_{x_3}(\delta\Phi^{\pm,n-1})(U_{a,7}^\pm+V_7^{\pm, n-1}),
%(?)
\end{cases}
\end{equation}
which implies
%\stepcounter{nonl}
\begin{equation}\label{f8.5}
\|({\cal E}_1^{\pm, n},{\cal E}_2^{\pm, n})\|_{s_0,T} \le
C\delta\theta_n^{s_0-1-\alpha}
\end{equation}
by using $\Delta_n=O(\theta_n^{-1})$, the inductive assumption for
(\ref{f7.2}), and Lemma 7.1. Thus we deduce
%\stepcounter{nonl}
\begin{equation}\label{f8.6}
\|S_{\theta_n}({\cal E}_1^{\pm, n},{\cal E}_2^{\pm, n})\|_{s,T} \le
C\delta\theta_n^{s-\alpha-1}\qquad\mbox{for any}\,\, s\ge s_0.
\end{equation}

The discussion for the commutators in (\ref{f8.3}) follows an
argument from \cite{cou2}. We now analyze the third term of
$V_2^{\pm, n+\frac 12}-S_{\theta_n}V_2^{\pm, n}$ given by
(\ref{f8.3}) in detail.

When $s\in [\alpha+1,s_1+2]$, we have
$$
\begin{array}{ll}
&\|\partial_{x_2}(\Psi_a^\pm
 +S_{\theta_n}\Phi^{\pm,n})S_{\theta_n}V_3^{\pm, n}\|_{s,T}\\[3mm]
& \le C\big(\|\partial_{x_2}(\Psi_a^\pm
+S_{\theta_n}\Phi^{\pm,n})\|_{L^\infty}\|S_{\theta_n}V_3^{\pm,
n}\|_{s,T} +\|\dot{\Psi}_a^\pm +S_{\theta_n}\Phi^{\pm,n}\|_{s+1,
T}\|S_{\theta_n}V_3^{\pm,n}\|_{L^\infty}\big)
 \le C\delta^2\theta_n^{s+1-\alpha},
\end{array}
$$
and
$$\begin{array}{ll}
&\|S_{\theta_n}(\partial_{x_2}(\Psi_a^\pm
 +\Phi^{\pm,n})V_3^{\pm, n})\|_{s,T}\\[3mm]
&\le C\theta_n^{s-\alpha}\|\partial_{x_2}(\Psi_a^\pm
+\Phi^{\pm,n})\|_{L^\infty}\|V_3^{\pm, n}\|_{\alpha,T}
+\|\dot{\Psi}_a^\pm+\Phi^{\pm,n}\|_{\alpha+1, T}\|V_3^{\pm,
n}\|_{L^\infty}\le C\delta^2\theta_n^{s+1-\alpha}
\end{array}
$$
by using (\ref{f7.1}), the induction assumption for (\ref{f7.2}),
and Lemma 7.1, which implies
\begin{equation}\label{f8.7}
\|\partial_{x_2}(\Psi_a^\pm+S_{\theta_n}\Phi^{\pm,n})S_{\theta_n}V_3^{\pm, n}-
S_{\theta_n}(\partial_{x_2}(\Psi_a^\pm+\Phi^{\pm,n})V_3^{\pm, n})\|_{s,T}
\le C\delta^2\theta_n^{s-\alpha+1}.
\end{equation}

When $s\in [s_0, \alpha]$, we have
\begin{equation}\label{f8.8}
\begin{array}{ll}
\|\partial_{x_2}(\Psi_a^\pm
\hspace{-.1in} &+S_{\theta_n}\Phi^{\pm,n})S_{\theta_n}V_3^{\pm, n}-
S_{\theta_n}\big(\partial_{x_2}(\Psi_a^\pm+\Phi^{\pm,n})V_3^{\pm, n}\big)\|_{s,T}\\[3mm]
& \le \|(I-S_{\theta_n})(\partial_{x_2}(\Psi_a^\pm+\Phi^{\pm,n})V_3^{\pm, n})\|_{s,T}
+
\|\partial_{x_2}(\Psi_a^\pm+\Phi^{\pm,n})(S_{\theta_n}-I)V_3^{\pm, n}\|_{s,T}\\[3mm]
& \hspace{.1in}
+\|\partial_{x_2}((S_{\theta_n}-I)\Phi^{\pm,n})S_{\theta_n}V_3^{\pm, n}\|_{s,T}\\[3mm]
& \le
C\theta_n^{s-\alpha}\big(\|\partial_{x_2}(\Psi_a^\pm+\Phi^{\pm,n})V_3^{\pm,
n}\|_{\alpha,T}
+\|\partial_{x_2}(\Psi_a^\pm+\Phi^{\pm,n})\|_{L^\infty}\|V_3^{\pm, n}\|_{\alpha,T}\big)\\[3mm]
& \hspace{.1in}
+C\theta_n^{3-\alpha}\|V_3^{\pm, n}\|_{\alpha,T}\|\Psi_a^\pm+\Phi^{\pm,n}\|_{s+1,T}
+\|\partial_{x_2}((S_{\theta_n}-I)\Phi^{\pm,n})S_{\theta_n}V_3^{\pm, n}\|_{s,T}\\[3mm]
& \le C\delta^2\theta_n^{s+1-\alpha}.\end{array}
\end{equation}

Together (\ref{f8.7}) with (\ref{f8.8}), it follows that
\begin{equation}\label{f8.9}
\|\partial_{x_2}(\Psi^{\pm}_a+S_{\theta_n}\Phi^{\pm,
n})S_{\theta_n}V_3^{\pm, n}-
S_{\theta_n}(\partial_{x_2}(\Psi^{\pm}_a+\Phi^{\pm, n})V_3^{\pm,
n})\|_{s,T} \le C\delta^2\theta_n^{s-\alpha+1} \quad\mbox{for
all}\,\, s\in [s_0, s_1+2].
\end{equation}
Other commutators appeared in (\ref{f8.3}) satisfy the estimates
similar to the above. Thus, we complete the proof.
%\P

\vspace{.1in} To show the error estimates given in Lemma 7.3, we
first show several lemmas dealing with different type errors.

\vspace{.05in}
{\bf Lemma 8.1.} {\it Let $\alpha\ge s_0+1\ge 6$. For the quadratic errors, we
have
%\stepcounter{nonl}
\begin{equation}\label{f8.10}
\begin{cases} \|e_{\pm,k}^{(1)}\|_{s,T}\le
C\delta^2\theta_k^{L_1(s)}\Delta_k  \qquad\mbox{for}\,\, s\in
[s_0-2, s_1-2],\cr \|\bar{e}_{\pm,k}^{(1)}\|_ {s,T}\le
C\delta^2\theta_k^{s-2\alpha+s_0-2}\Delta_k\qquad\mbox{for}\,\, s\in
[s_0-1, s_1-1],\cr \|\tilde{e}_{k}^{(1)}\|_{H^s(b\Omega_T)}\le
C\delta^2\theta_k^{s-2\alpha+s_0-1}\Delta_k  \qquad\mbox{for}\,\,
s\in [s_0-2, s_1-2]
\end{cases}
\end{equation}
for all $k\le n-1$, where
%\stepcounter{nonl}
\begin{equation}\label{f8.11}
L_1(s)=\max((s+2-\alpha)_++2(s_0-\alpha)-3, s+s_0-1-2\alpha, s+2s_0-3-3\alpha, s_0-\alpha-3).
\end{equation}}

%\vspace{.1in}
{\it Proof.} The proof is divided into three steps.

{\it Step 1}. From the definition of $e_{\pm,k}^{(1)}$, we first
have %\stepcounter{nonl}
\begin{equation}\label{f8.12}
e_{\pm,k}^{(1)}=\int_0^1(1-\tau)L{''}_{(U_a^\pm+V^{\pm,
k}+\tau\delta V^{\pm, k}; \Psi_a^\pm+\Phi^{\pm, k}+\tau\delta
\Phi^{\pm, k})}\big((\delta V^{\pm, k},\delta \Phi^{\pm, k}),
(\delta V^{\pm, k},\delta \Phi^{\pm, k})\big)d\tau.
\end{equation}

{}From (\ref{f7.1}) and Lemma 7.1, we find
\begin{equation}\label{f8.13}
\sup_{0\le \tau\le 1}(\|\dot{U}_a^\pm+V^{\pm,k}+\tau\delta
V^{\pm,k}\|_ {W^{1,\infty}(\Omega_T)} +
\|\dot{\Psi}_a^\pm+\Phi^{\pm,k}+\tau\delta \Phi^{\pm,k}\|_
{W^{1,\infty}(\Omega_T)})\le C\delta.
\end{equation}

On the other hand, we
%obviously
have
%\stepcounter{nonl}
\begin{equation}\label{f8.14}
\begin{aligned}
&\|L''_{(U^\pm, \Psi^\pm)}((V^{\pm,1}, \Phi^{\pm,1}),(V^{\pm,2},
\Phi^{\pm,2}))\|_{s,T}\\
& \le C\Big(\|(\dot{U}^\pm, \dot{\Psi}^\pm)\|_{s+2,T}\|(V^{\pm,1},
\Phi^{\pm,1})\|_{W^{1,\infty}} \|(V^{\pm,2},\Phi^{\pm,2})\|_{W^{1,\infty}}\\
&\qquad\,+ \|(V^{\pm,1}, \Phi^{\pm,1})\|_{s+2,T} \|(V^{\pm,2},
\Phi^{\pm,2})\|_{W^{1,\infty}}+ \|(V^{\pm,1},
\Phi^{\pm,1})\|_{W^{1,\infty}} \|(V^{\pm,2},
\Phi^{\pm,2})\|_{s+2,T}\Big).
\end{aligned}
\end{equation}
Therefore, we obtain
$$
\begin{aligned}
\|e^{(1)}_{\pm, k}\|_{s,T} \le
&C\Big(\delta^2\theta_k^{2(s_0-1-\alpha)}\Delta_k^2
(\delta+\|(V^{\pm,k}, \Phi^{\pm,k})\|_{s+2,T}+ \|(\delta V^{\pm,k},
\delta
\Phi^{\pm,k})\|_{s+2,T})\\
&\quad\,\,+ \delta\theta_k^{s_0-1-\alpha} \Delta_k\|(\delta
V^{\pm,k}, \delta\Phi^{\pm,k})\|_{s+2,T}\Big),
\end{aligned}
$$
which implies
%\stepcounter{nonl}
\begin{equation}\label{f8.15}
\begin{cases} \|e^{(1)}_{\pm,k}\|_{s,T} \le
C\delta^2\theta_k^{L(s)}\Delta_k \qquad {\rm when}~s+2\neq \alpha,
~s\le s_1-2,\cr \|e^{(1)}_{\pm,k}\|_{s,T} \le
C\delta^2\theta_k^{\max(s_0-\alpha-3, 2s_0-2-2\alpha)}\Delta_k \qquad {\rm when}~s+2=\alpha,
\end{cases}
\end{equation}
by Lemma 7.1, the inductive assumption for (\ref{f7.2}), and
$\Delta_k=O({\theta_k}^{-1})$, where
$$
L(s)=\max((s+2-\alpha)_++2(s_0-\alpha)-3, s+s_0-1-2\alpha,
s+2s_0-3-3\alpha).
$$
Thus, we conclude the first result in (\ref{f8.10}).

{\it Step 2}.  Obviously, from the definition  of $\bar{e}_{\pm,
k}^{(1)}$ in (\ref{f6.29}), we have
%\stepcounter{nonl}
\begin{equation}\label{f8.16}
\|\bar{e}_{\pm,k}^{(1)}\|_{s,T}\le C\left(\|\delta
V^{\pm,k}\|_{s,T}\|\delta\Phi^{\pm,k}\|_{W^{1,\infty}}+ \|\delta
V^{\pm,k}\|_{L^\infty}\|\delta\Phi^{\pm,k}\|_{s+1, T}\right),
\end{equation}
which implies the second result in (\ref{f8.10}).

{\it Step 3}. Since $\bar{e}_{\pm,
k}^{(1)}|_{x_1=0}=(\tilde{e}_k)_1^\pm$, we obtain that
$(\tilde{e}_k)_1^\pm$ satisfy the estimate given in (\ref{f8.10}).
Moreover, from the definition of $\tilde{e}_k^{(1)}$, we have
%\stepcounter{nonl}
\begin{equation}\label{f8.17}
(\tilde{e}_{k}^{(1)})_2^\pm=
-\partial_{x_2}(\delta\phi^k)\delta V_6^{\pm,k}-
\partial_{x_3}(\delta\phi^k)\delta V_7^{\pm,k}, \quad
(\tilde{e}_{k}^{(1)})_3=\frac 12(|\delta V_H^{+,k}|^2-|\delta
V_H^{-,k}|^2),
\end{equation}
which implies
%\stepcounter{nonl}
\begin{equation}\label{f8.18}
\begin{cases}
\|(\tilde{e}_{k}^{(1)})_2^\pm\|_{H^s(b\Omega_T)}\le C\big(\|\delta
\phi^{k}\|_{H^{s+1}(b\Omega_T)}\|\delta V^{\pm, k}\|_{L^\infty}+
\|\delta \phi^{k}\|_{W^{1,\infty}(b\Omega_T)}\|\delta V^{\pm,
k}\|_{H^s(b\Omega_T)}\big),\cr
\|(\tilde{e}_{k}^{(1)})_3\|_{H^s(b\Omega_T)}\le C\|\delta
V^{\pm,k}\|_{L^\infty}\|\delta V^{\pm, k}\|_{H^s(b\Omega_T)},
\end{cases}
\end{equation}
yielding the last estimate of (\ref{f8.10})
%immediately
by the inductive assumption for (\ref{f7.2}).
%\P

\vspace{.1in}
For the first substitution errors, we have

\vspace{.1in} {\bf Lemma 8.2.} {\it Let $\alpha\ge s_0+1\ge 6$.
Then, for all $k\le n-1$, we have
%\stepcounter{nonl}
\begin{equation}\label{f8.19}
\begin{cases} \|e_{\pm,k}^{(2)}\|_{s,T}\le
C\delta^2\theta_k^{L_2(s)}\Delta_k \qquad\mbox{for}\,\, s\in
[s_0-2,s_1-2],\cr \|\bar{e}_{\pm,k}^{(2)}\|_ {s,T}\le
C\delta^2\theta_k^{s-2\alpha+s_0}\Delta_k \qquad\mbox{for}\,\, s\in
[s_0-1, s_1-1],\cr \|\tilde{e}_{k}^{(2)}\|_{H^s(b\Omega_T)}\le
C\delta^2\theta_k^{s-2\alpha+s_0+1}\Delta_k \qquad\mbox{for}\,\,
s\in [s_0-2, s_1-2],
\end{cases}
\end{equation}
where
%\stepcounter{nonl}
\begin{equation}\label{f8.20}
L_2(s)=\begin{cases} (s+2-\alpha)_++s_0-\alpha-1
\qquad\mbox{for}\,\, \alpha\neq s+2,\cr \max(5-\alpha, s_0-1-\alpha)
\qquad \mbox{for}\,\, \alpha=s+2.
\end{cases}
\end{equation}}

{\it Proof.} The proof is divided into three steps.

{\it Step 1}. From the definition of $e_{\pm,k}^{(2)}$, we
have
%\stepcounter{nonl}
\begin{equation}\label{f8.21}
\begin{aligned}
e_{\pm,k}^{(2)}=&\int_0^TL{''}_{(U_a^\pm+S_{\theta_k}V^{\pm,
k}+\tau(1-S_{\theta_k})V^{\pm, k}; \Psi_a^\pm+S_{\theta_k}\Phi^{\pm,
k}+\tau(1-S_{\theta_k})\Phi^{\pm, k})}\\
&\qquad \times ((\delta V^{\pm, k},\delta \Phi^{\pm,
k}),((1-S_{\theta_k})V^{\pm, k},(1-S_{\theta_k})\Phi^{\pm,
k}))d\tau.
\end{aligned}
\end{equation}

As in (\ref{f8.13}), from the inductive assumption for (\ref{f7.2}),
we find
\begin{equation}\label{f8.22}
\sup_{0\le \tau\le 1}\|\dot{U}_a^\pm+S_{\theta_k}V^{\pm,k}
+\tau(1-S_{\theta_k})V^{\pm,k}\|_{W^{1,\infty}(\Omega_T)}\le
C\delta,
\end{equation}
and
%\stepcounter{nonl}
\begin{equation}\label{f8.23}
\sup_{0\le \tau\le 1}\|\dot{\Psi}_a^\pm
+S_{\theta_k}\Phi^{\pm,k}+\tau(1-S_{\theta_k})\Phi^{\pm,k}\|_
{W^{1,\infty}(\Omega_T)}\le C\delta.
\end{equation}
Then we use (\ref{f8.14}) in (\ref{f8.21}) to obtain
%\stepcounter{nonl}
\begin{equation}\label{f8.24}
\begin{aligned}
\|e_{\pm,k}^{(2)}\|_{s,T} \le& C\big(\|(\delta V^{\pm, k},\delta
\Phi^{\pm, k})\|_{W^{1,\infty}} \|(1-S_{\theta_k})( V^{\pm,
k},\Phi^{\pm, k})\|_{W^{1,\infty}} \cr &\quad
\times(\|(\dot{U}_a^\pm, \dot{\Psi}^\pm_a)\|_{s+2, T}
+\|S_{\theta_k}(V^{\pm, k},\Phi^{\pm, k})\|_{s+2,T}+
\|(1-S_{\theta_k})(V^{\pm, k},\Phi^{\pm, k})\|_{s+2,T})\cr &\quad
+\|(\delta V^{\pm, k},\delta \Phi^{\pm, k})\|_{s+2,T}
\|(1-S_{\theta_k})(V^{\pm, k},\Phi^{\pm, k})\|_{W^{1,\infty}}\cr
&\quad +\|(\delta V^{\pm, k},\delta \Phi^{\pm, k})\|_{W^{1,\infty}}
\|(1-S_{\theta_k})(V^{\pm, k},\Phi^{\pm, k})\|_{s+2, T}\big).
\end{aligned}
\end{equation}

Using the properties of the smoothing operators, the inductive assumption
for (\ref{f7.2}), and Lemma 7.1 in (\ref{f8.24}), we conclude the first
estimate in (\ref{f8.19}) when $s\le s_1-2$.

\medskip
{\it Step 2.} From the definition of $\bar{e}_{\pm, k}^{(2)}$ in
(\ref{f6.30}), we have
\begin{equation}\label{f8.25}
\begin{aligned}\|\bar{e}_{\pm,k}^{(2)}\|_{s,T}\le& C\Big(\|\delta V^{\pm,
k}\|_{s,T}\|\nabla_{(x_2,x_3)}(1-S_{\theta_k})\Phi^{\pm,
k}\|_{L^{\infty}}+ \|(1-S_{\theta_k})\Phi^{\pm, k}\|_{s+1,
T}\|\delta V^{\pm, k}\|_{L^\infty} \cr &\quad +
\|(1-S_{\theta_k})V^{\pm, k}\|_{s, T}\|\nabla_{(x_2,x_3)}\delta
\Phi^{\pm, k}\|_{L^{\infty}}+ \|\delta\Phi^{\pm, k}\|_{s+1,
T}\|(1-S_{\theta_k})V^{\pm, k}\|_{L^\infty}\Big),
\end{aligned}
\end{equation}
which implies the second result in (\ref{f8.19}) when $s\le s_1-1$
by the inductive assumption for (\ref{f7.2}) and Lemma 7.1 in
(\ref{f8.25}).

\medskip
{\it Step 3.} Noting that
$$
(\tilde{e}^{(2)}_k)_1^\pm=\bar{e}_{\pm, k}^{(2)}|_{x_1=0}
$$
from the above discussion, we conclude that
$(\tilde{e}^{(2)}_k)_1^\pm$ satisfy the estimate given in (\ref{f8.19})
when $s\le s_1-2$.

{}From the definition of $\tilde{e}_k^{(2)}$, we have
%\stepcounter{nonl}
\begin{equation}\label{f8.26}
\begin{cases}
(\tilde{e}_{k}^{(2)})_2^\pm
=\partial_{x_2}((S_{\theta_k}-I)\phi^k)\delta V_6^{\pm, k}
+\partial_{x_3}((S_{\theta_k}-I)\phi^k)\delta V_7^{\pm, k}
+(S_{\theta_k}-I)V_6^{\pm, k}\partial_{x_2}(\delta\phi^k)\cr
\qquad\qquad\, +(S_{\theta_k}-I)V_7^{\pm,
k}\partial_{x_3}(\delta\phi^k),\cr (\tilde{e}_{k}^{(2)})_3
=(I-S_{\theta_k})V_H^{+, k}\delta V_H^{+, k}-
(I-S_{\theta_k})V_H^{-, k}\delta V_H^{-, k},
\end{cases}
\end{equation}
which implies
%\stepcounter{nonl}
\begin{equation}\label{f8.27}
\begin{cases} \|(\tilde{e}_{k}^{(2)})_2^\pm\|_{H^s(b\Omega_T)}\le
C\big(\|\delta V^{\pm,
k}\|_{H^s(b\Omega_T)}\|(1-S_{\theta_k})\phi^{k}\|_{W^{1,\infty}}+
\|(1-S_{\theta_k})\phi^{k}\|_{H^{s+1}(b\Omega_T)}\|\delta V^{\pm,
k}\|_{L^\infty}\cr \hspace{1.4in}+ \|(1-S_{\theta_k})V^{\pm,
k}\|_{H^s(b\Omega_T)}\|\delta \phi^{ k}\|_{W^{1,\infty}}+
\|\delta\phi^{k}\|_{H^{s+1}(b\Omega_T)}\|(1-S_{\theta_k})V^{\pm,
k}\|_{L^\infty}\big),\cr
\|(\tilde{e}_{k}^{(2)})_3\|_{H^s(b\Omega_T)}\le C\big(\|\delta
V^{\pm,
k}\|_{H^s(b\Omega_T)}\|(1-S_{\theta_k})V^{\pm,k}\|_{L^{\infty}}+
\|(1-S_{\theta_k})V^{\pm,k}\|_{H^{s}(b\Omega_T)}\|\delta V^{\pm,
k}\|_{L^\infty}\big).
\end{cases}
\end{equation}

Using the inductive assumption for (\ref{f7.2}) and Lemma 7.1 in
(\ref{f8.27}), we conclude the last result in (\ref{f8.19}) when
$s\le s_1-2$.
 %\P

\vspace{.1in} The representations of the second substitution errors
$(e^{(3)}_{\pm, k},
\tilde{e}^{(3)}_{k})$ are similar to those of
$(e^{(2)}_{\pm, k}, \tilde{e}^{(2)}_{k})$. Thus, using Lemma 7.2 and the same argument
as the proof of Lemma 8.2, we conclude

\vspace{.1in} {\bf Lemma 8.3.} {\it Let $\alpha\ge s_0+1\ge 6$. For the second substitution
errors, we have
%\stepcounter{nonl}
\begin{equation}\label{f8.28}
\begin{cases} \|e_{\pm,k}^{(3)}\|_{s,T}\le
C\delta^2\theta_k^{L_3(s)}\Delta_k \qquad\mbox{for}\,\,s\in [s_0-2,
s_1-2],\cr \|\tilde{e}_{k}^{(3)}\|_{H^s(b\Omega_T)}\le
C\delta^2\theta_k^{s+s_0+2-2\alpha}\Delta_k \qquad\mbox{for}\,\,
s\in [s_0-2, s_1-2]
\end{cases}
\end{equation}
for all $k\le n-1$, where
$L_3(s)=\max(s+s_0+2-2\alpha, s+2s_0+3-3\alpha,
(s+2-\alpha)_++2(s_0-\alpha)+1).$}
%\end{equation}}

{\bf Lemma 8.4.} {\it Let $s_0\ge 5$ and $\alpha\ge s_0+3$. For the last errors, we have
%\stepcounter{nonl}
\begin{equation}\label{f8.30}
\begin{cases} \|e_{\pm,k}^{(4)}\|_{s,T}\le
C\delta^2\theta_k^{L_4(s)}\Delta_k \qquad\mbox{for}\,\, s\in [s_0,
s_1-4],\cr \|\bar{e}_{\pm,k}^{(4)}\|_ {s,T}\le
C\delta^2\theta_k^{L_5(s-1)}\Delta_k \qquad\mbox{for}\,\, s\in [s_0,
s_1-3],\cr \|\tilde{e}_{k}^{(4)}\|_{H^s(b\Omega_T)}\le
C\delta^2\theta_k^{L_5(s)}\Delta_k \qquad\mbox{for}\,\, s\in [s_0-1,
s_1-4]
\end{cases}
\end{equation}
for all $k\le n-1$, where
$$
L_4(s)=
\begin{cases}
\max((s+2-\alpha)_++2+2(s_0-\alpha), s+4+s_0-2\alpha) \qquad {\rm
if}\quad \alpha\neq s+2, s+4,\cr \max(s_0-\alpha, 2(s_0-\alpha)+4)
\qquad {\rm if}\quad \alpha=s+4,\cr s_0+2-\alpha) \qquad {\rm
if}\quad \alpha=s+2,
\end{cases}
$$
and
$$
L_5(s)=
\begin{cases}
\max((s+3-\alpha)_++2(s_0-\alpha), s+3+s_0-2\alpha) \qquad {\rm
if}\quad \alpha\neq s+3, s+4,\cr s_0-\alpha \qquad {\rm if}\quad
\alpha=s+3,\cr s_0-1-\alpha) \qquad {\rm if}\quad \alpha=s+4.
\end{cases}
$$}

{\it Proof.} The proof is divided into three steps.

{\it Step 1}. Set
\begin{equation}\label{f8.32}
R^\pm_k=\partial_{x_1}\big(L(U_a^\pm+V^{\pm, k+\frac 12},\Psi_a^\pm
+S_{\theta_k}\Phi^{\pm, k})(U_a^\pm+V^{\pm, k+\frac 12})\big).
\end{equation}
Then we have
$$
\begin{aligned}
\|R^\pm_k\|_{s,T}\le& \|L(U_a^\pm+V^{\pm, k+\frac
12},\Psi_a^\pm+S_{\theta_k}\Phi^{\pm, k})(U_a^\pm+V^{\pm, k+\frac
12})\cr &\,\, - L(U_a^\pm+V^{\pm, k},\Psi_a^\pm+\Phi^{\pm,
k})(U_a^\pm+V^{\pm, k})\|_{s+2,T}\cr &\,\, + \|{\cal L}(V^{\pm,
k},\Phi^{\pm, k})V^{\pm, k}-f_a^\pm\|_{s+2,T},
\end{aligned}
$$
which implies that, for all $s\in [s_0-2, s_1-4]$,
%\stepcounter{nonl}
\begin{equation}\label{f8.33}
\begin{aligned}
\|R^\pm_k\|_{s,T}\le & C\Big(\|V^{\pm, k+\frac 12}-V^{\pm,
k}\|_{W^{1,\infty}}(\|\dot{U}_a^\pm+V^{\pm, k+\frac 12}\|_{s+2, T}
+\|\dot{\Psi}_a^\pm+S_{\theta_k}\Phi^{\pm, k}\|_{s+4, T})\\
&\quad + \|\dot{U}_a^\pm+V^{\pm, k}\|_{s+4, T}(\|V^{\pm, k+\frac
12}-V^{\pm,
k}\|_{L^\infty}+\|(I-S_{\theta_k})\Phi^{\pm, k}\|_{W^{1,\infty}})\\
&\quad + \|\dot{U}_a^\pm+V^{\pm, k}\|_{W^{1,\infty}}(\|V^{\pm,
k+\frac 12}-V^{\pm,
k}\|_{s+2,T}+\|(I-S_{\theta_k})\Phi^{\pm, k}\|_{s+4, T})\\
&\quad+\|V^{\pm, k+\frac 12}-V^{\pm, k}\|_{s+4,T}+ \|{\cal
L}(V^{\pm, k},\Phi^{\pm, k})V^{\pm, k}-f_a^\pm\|_{s+2,T}\Big)\cr \le
&
\begin{cases}
C\delta^2\theta_k^{(s+4-\alpha)_++s_0+1-\alpha}+C\delta\theta_k^{s+5-\alpha}
 \quad\mbox{for}\,\, s+4\neq \alpha,\cr
C\delta^2\theta_k^{s_0+2-\alpha}+C\delta\theta_k \quad
\mbox{for}\,\, s+4=\alpha.
\end{cases}
 \end{aligned}
\end{equation}
Thus we find that
$$
e_{\pm, k}^{(4)}
=\frac{R_k^\pm \delta\Phi^{\pm,k}}{\partial_{x_1}(\Psi_a^\pm
+S_{\theta_k}\Phi^{\pm,k})}
$$
satisfy
%\stepcounter{nonl}
\begin{equation}\label{f8.34}
\|e_{\pm, k}^{(4)}\|_{s,T}\le
C\big(\|R^\pm_k\|_{s_0-2,T}(\delta\theta_k^{s-1-\alpha}\Delta_k
+\delta\theta_k^{s_0-1-\alpha}\Delta_k
(\delta+\delta\theta_k^{(s+2-\alpha)_+}))
+\delta\theta_k^{s_0-1-\alpha}\Delta_k\|R_k^\pm\|_{s,T}\big),
\end{equation}
when $s+2\neq \alpha$, and
\begin{equation}\label{f8.35}
\|e_{\pm, k}^{(4)}\|_{s,T}\le
C\big(\|R^\pm_k\|_{s_0-2,T}(\delta\theta_k^{s-1-\alpha}\Delta_k
+\delta\theta_k^{s_0-1-\alpha}\Delta_k(\delta+\delta\log\theta_k))
+\delta\theta_k^{s_0-1-\alpha}\Delta_k\|R_k^\pm\|_{s,T}\big),
\end{equation}
when $s+2=\alpha$, which yields the first estimate in (\ref{f8.30})
for any $s\in [s_0,s_1-4]$, provided $\alpha\ge s_0+3$ by using
(\ref{f8.33}).

\medskip
{\it Step 2}. Set
%\stepcounter{nonl}
\begin{equation}\label{f8.36}
R^b_k=B(U_a^\pm+V^{\pm, k+\frac 12},\psi_a+S_{\theta_k}\phi^{k}).
\end{equation}
Then we have
$$
%\begin{aligned}
\|R^b_k\|_{H^s(b\Omega_T)}\le \|B(U_a^\pm+V^{\pm, k+\frac
12},\psi_a+S_{\theta_k}\phi^{k})-B(U_a^\pm+V^{\pm,
k},\psi_a+\phi^{k})\|_{H^s(b\Omega_T)}
%\cr \hspace{1.0in}
+\|{\cal
B}(V^{\pm, k},\phi^{k})\|_{H^s(b\Omega_T)},
%\end{aligned}
$$
which
implies
%\stepcounter{nonl}
\begin{equation}\label{f8.37}
\begin{aligned}
&\|(R^b_k)_1^\pm\|_{H^s(b\Omega_T)}\\
&\le C\Big( \|(S_{\theta_k}-1)\phi^{k}\|_{H^{s+1}(b\Omega_T)}+
\|V_2^{\pm,k+\frac
12}-V_2^{\pm,k}\|_{H^{s}(b\Omega_T)}+\|(S_{\theta_k}-I)V^{\pm,k}\|_{L^\infty}
\|\dot{\psi}_a+S_{\theta_k}\phi^{k}\|_{H^{s+1}(b\Omega_T)} \cr
%\hspace{.3in}
&\qquad\,\,+\|(S_{\theta_k}-I)V^{\pm,k}\|_{H^{s}(b\Omega_T)}
\|\dot{\psi}_a+S_{\theta_k}\phi^{k}\|_{W^{1,\infty}}
+\|(S_{\theta_k}-1)\phi^{k}\|_{H^{s+1}(b\Omega_T)}
\|\dot{U}_a^\pm+V^{\pm,k}\|_{L^{\infty}} \cr
%\hspace{.3in}
&\qquad\,\,+\|(S_{\theta_k}-1)\phi^{k}\|_{W^{1,\infty}(b\Omega_T)}
\|\dot{U}_a^\pm+V^{\pm,k}\|_{H^s(b\Omega_T)}
+\|{\cal B}(V^{\pm, k},\phi^{k})\|_{H^s(b\Omega_T)}\Big).
\end{aligned}
\end{equation}
Therefore, we have the estimate
%\stepcounter{nonl}
\begin{equation}\label{f8.38}
\|(R^b_k)_1^\pm\|_{H^s(b\Omega_T)}\le
\begin{cases}
 C\delta\theta_k^{\max((s+2-\alpha)_++s_0-\alpha, s+2-\alpha)}
\qquad\mbox{if}\,\, \alpha\neq s+1, s+2,\cr
 C\delta\theta_k
\qquad\mbox{if}\,\, \alpha=s+1,\cr
 C\delta
\qquad\mbox{if}\,\, \alpha=s+2\end{cases}
\end{equation}
for all $s\in [s_0-1, s_1-2]$. Thus we find that
$$
(\tilde{e}_{k}^{(4)})_1^\pm
=-\frac{\partial_{x_1}(R_k^b)_1^\pm}{\partial_{x_1}(\Psi_a^\pm
+S_{\theta_k}\Phi^{\pm,k})_{x_1=0}} \delta\phi^{k}
$$
satisfy
%\stepcounter{nonl}
\begin{equation}\label{f8.39}
\begin{array}{ll}
\|(\tilde{e}_{k}^{(4)})_1^\pm\|_{H^s(b\Omega_T)}\le
&C\big(\|(R^b_k)_1^\pm\|_{H^4(b\Omega_T)}(\delta\theta_k^{s-\alpha}\Delta_k
+\delta^2\theta_k^{s_0-1-\alpha}\Delta_k
\theta_k^{(s+3-\alpha)_+})\\
&\quad\,\,
+\delta^2\theta_k^{s_0-1-\alpha}\Delta_k\|(R^b_k)_1^\pm\|_{H^{s+2}(b\Omega_T)}\big),
\end{array}
\end{equation}
when $s+3\neq \alpha$, and
\begin{equation}\label{f8.40}
\|(\tilde{e}_{k}^{(4)})_1^\pm\|_{H^s(b\Omega_T)}\le
C\big(\|(R^b_k)_1^\pm\|_{H^4(b\Omega_T)}(\delta\theta_k^{s-\alpha}\Delta_k
+\delta^2\theta_k^{s_0-1-\alpha}\Delta_k \log\theta_k)
+\delta^2\theta_k^{s_0-1-\alpha}\Delta_k\|(R^b_k)_1^\pm\|_{H^{s+2}(b\Omega_T)}\big),
\end{equation}
when $s+3=\alpha$, which yields the estimate of
$(\tilde{e}_{k}^{(4)})_1^\pm$ given in (\ref{f8.30}) for any $s\in
[s_0-1, s_1-4]$ by using (\ref{f8.38}).

Estimate (\ref{f8.30}) for the other components of
$\tilde{e}_{k}^{(4)}$ can be proved by the same argument as above.

\medskip
{\it Step 3.} Noting that the restriction of $\bar{e}_{\pm,k}^{(4)}$
on $\{x_1=0\}$ is the same as $(\tilde{e}_{k}^{(4)})_1^\pm$, we
conclude the estimate of $\bar{e}_{\pm,k}^{(4)}$ in (\ref{f8.30})
for any $s\in [s_0, s_1-3]$.
%\P

\vspace{.1in} {\sc Proof of Lemma 7.3:} Summarizing all the results
in Lemmas 8.1--8.4, we obtain the estimates given in Lemma 7.3.

\vspace{.1in} {\sc Proof of Lemma 7.5.} {}From the definitions of
$(f_n^\pm, g_n, h_n^\pm)$ given in (\ref{f6.22}) and (\ref{f6.35}),
respectively, we have
%\stepcounter{nonl}
\begin{equation}\label{f8.41}
\begin{cases}
f_n^\pm=(S_{\theta_n}-S_{\theta_{n-1}})f_a^\pm-(S_{\theta_n}-S_{\theta_{n-1}})E_{\pm,
n-1} -S_{\theta_{n}}e_{\pm, n-1},\cr
g_n=-(S_{\theta_n}-S_{\theta_{n-1}})\tilde{E}_{n-1}-S_{\theta_{n}}\tilde{e}_{n-1},\cr
h_n^\pm=(S_{\theta_{n-1}}-S_{\theta_{n}})\bar{E}_{\pm,n-1}
-S_{\theta_{n}}\bar{e}_{\pm,n-1}. \end{cases}
\end{equation}

Using the properties of the smoothing operators, we have
\begin{equation}\label{f8.42}
\begin{cases}
\|(S_{\theta_n}-S_{\theta_{n-1}})f_a^\pm\|_{s,T} \le
C\theta_n^{s-\tilde{s}-1}\Delta_n\|f_a^\pm\|_{\tilde{s},T}
\qquad\mbox{for}\,\, \tilde{s}\ge 0,\cr
\|(S_{\theta_n}-S_{\theta_{n-1}})E_{\pm, n-1}\|_{s,T}\le
C\theta_n^{s-\tilde{s}-1}\Delta_n\|E_{\pm, n-1}\|_{\tilde{s},T}\le
C\delta^2\theta_n^{s-\tilde{s}}\Delta_n \qquad\mbox{for}\,\,
\tilde{s}\in [s_0, s_1-4],\cr \|S_{\theta_{n}}{e}_{\pm,n-1}\|_{s,T}
\le C\theta_n^{(s-\tilde{s})_+}\|e_{\pm, n-1}\|_{\tilde{s},T}\le C
\delta^2\theta_n^{(s-\tilde{s})_++L^1(\tilde{s})}\Delta_n
 \qquad\mbox{for}\,\, \tilde{s}\in [s_0, s_1-4].
 \end{cases}
\end{equation}
Estimate $(\ref{f8.42})_3$ can be represented as the following three
cases.

{\it Case 1: $s=s_0$}. For this case, if we choose $\tilde{s}\in
[s_0, \alpha-5]$, then
$$
L^1(\tilde{s})=\max(s_0-\alpha-1,
\tilde{s}+s_0+4-2\alpha)=s_0-\alpha-1,
$$
which follows from
$(\ref{f8.42})_3$  that
\begin{equation}\label{f8.43}
\|S_{\theta_{n}}{e}_{\pm,n-1}\|_{s_0,T}\le C\delta^2
\theta_n^{s_0-\alpha-1}\Delta_n.
\end{equation}

{\it Case 2: $s=s_0+1$ or $s=s_0+2$}. For these two cases, if we
choose $\tilde{s}=\alpha-4\ge s_0+1$, then
$$
L^1(\tilde{s})=\max(s_0-\alpha, 2(s_0-\alpha)+4)=s_0-\alpha,
$$
which
implies
\begin{equation}\label{f8.44}
\|S_{\theta_{n}}{e}_{\pm,n-1}\|_{s,T}\le C\delta^2
\theta_n^{s-\alpha-1}\Delta_n.
\end{equation}

{\it Case 3: $s\ge s_0+3$}. For this case, if we choose
$\tilde{s}=\alpha-2\ge s_0$, then
$$(s-\tilde{s})_++L^1(\tilde{s})=(s+2-\alpha)_++s_0+2-\alpha\le s-\alpha-1$$
by using $\alpha\ge s_0+5$,
which follows from  $(\ref{f8.42})_3$  that
\begin{equation}\label{f8.45}
\|S_{\theta_{n}}{e}_{\pm,n-1}\|_{s,T}\le C\delta^2
\theta_n^{s-\alpha-1}\Delta_n.
\end{equation}

On the other hand, we have
\begin{equation}\label{f8.46}
\|(S_{\theta_{n}}-S_{\theta_{n-1}}){E}_{\pm,n-1}\|_{s,T}\le C\delta^2
\theta_n^{s-\alpha-1}\Delta_n
\end{equation}
by setting $\tilde{s}=\alpha+1$ in $(\ref{f8.42})_2$ if $s_1\ge
\alpha+5$.

\medskip
Together (\ref{f8.43})--(\ref{f8.46}) with (\ref{f8.41}), it follows
\begin{equation}\label{f8.47}
\|f_n^\pm\|_{s,T}\le
C\theta_n^{s-\alpha-1}\Delta_n(\delta^2+\|f_a^\pm\|_{\alpha, T})
\qquad\mbox{for all}\,\,s\ge s_0.
\end{equation}

Similarly, we have
\begin{equation}\label{f8.48}
\begin{cases}
\|(S_{\theta_n}-S_{\theta_{n-1}})\tilde{E}_{n-1}\|_{H^{s+1}(b\Omega_T)}\le
C\theta_n^{s-\tilde{s}}\Delta_n\|\tilde{E}_{n-1}\|_{H^{\tilde{s}}(b\Omega_T)}
\le C\delta^2\theta_n^{s+1-\tilde{s}}\Delta_n
%\qquad\qquad\qquad\qquad\qquad\qquad\qquad\qquad\qquad\qquad\qquad\qquad
\,\,\,\mbox{for}\,\,\tilde{s}\in [s_0-1, s_1-4],\cr
\|S_{\theta_{n}}\tilde{e}_{n-1}\|_{H^{s+1}(b\Omega_T)}\le
C\theta_n^{(s+1-\tilde{s})_+}\|\tilde{e}_{n-1}\|_{H^{\tilde
s}(b\Omega_T)}\le
C\theta_n^{(s+1-\tilde{s})_++L^3(\tilde{s})}\delta^2\Delta_n \,\,\,
\mbox{for}\,\, \tilde{s}\in [s_0-1, s_1-4].
\end{cases}
\end{equation}
Estimate $(\ref{f8.48})_2$ can be represented as the following three
cases.

{\it Case 1: $s=s_0$}. For this case, if we choose
$\tilde{s}=\alpha-4$, then $L^3(\alpha-4)=s_0-\alpha-1$, which
follows from  $(\ref{f8.48})_3$  that
%\begin{equation}\label{f8.49}
$\|S_{\theta_{n}}\tilde{e}_{n-1}\|_{H^{s_0+1}(b\Omega_T)}\le
C\delta^2 \theta_n^{s_0-\alpha-1}\Delta_n$;
%\end{equation}

{\it Case 2: $s_0+1\le s\le\alpha-3$}. For this case, if we choose
$\tilde{s}=\alpha-3$, then $L^3(\alpha-3)=s_0-\alpha$, which implies
%\begin{equation}\label{f8.50}
$\|S_{\theta_{n}}\tilde{e}_{n-1}\|_{H^{s+1}(b\Omega_T)}\le C\delta^2
\theta_n^{s_0-\alpha}\Delta_n \le C\delta^2
\theta_n^{s-\alpha-1}\Delta_n$;
%\end{equation}

{\it Case 3: $s\ge \alpha-2$}. For this case, if we choose
$\tilde{s}=\alpha-5\ge s_0$, then $L^3(\tilde{s})=s_0-2-\alpha$,
which implies
%\begin{equation}\label{f8.51}
$\|S_{\theta_{n}}\tilde{e}_{n-1}\|_{H^{s+1}(b\Omega_T)}\le C\delta^2
\theta_n^{(s-\tilde{s})++L^3(\tilde{s})}\Delta_n\le
C\delta^2\theta_n^{s-\alpha-1}\Delta_n.$
%\end{equation}

\medskip
In summary, we obtain
\begin{equation}\label{f8.52}
\|S_{\theta_{n}}\tilde{e}_{n-1}\|_{H^{s+1}(b\Omega_T)}\le
C\delta^2\theta_n^{s-\alpha-1}\Delta_n \qquad\mbox{for all}\,\, s\ge
s_0.
\end{equation}

{}From the assumption $s_1\ge \alpha+6$, it is possible to let
$\tilde{s}=\alpha+2$ in $(\ref{f8.48})_1$, which yields
\begin{equation}\label{f8.53}
\|g_n\|_{H^{s+1}(b\Omega_T)}\le
C\delta^2\theta_n^{s-\alpha-1}\Delta_n \qquad\mbox{for all}\,\, s\ge
s_0,
\end{equation}
by using $(\ref{f8.48})_1$ and $(\ref{f8.52})$ in $(\ref{f8.41})_2$.

Furthermore, we have
\begin{equation}\label{f8.54}
\begin{cases}
\|(S_{\theta_n}-S_{\theta_{n-1}})\bar{E}_{\pm,n-1}\|_{s,T}\le
C\theta_n^{s-\tilde{s}-1}\Delta_n\|\bar{E}_{\pm,n-1}\|_{\tilde{s},T}
\le C\delta^2\theta_n^{s-\tilde{s}-1}\Delta_n
%\cr\qquad\qquad\qquad\qquad\qquad\qquad\qquad\qquad\qquad\qquad\qquad\qquad
\qquad\mbox{for}\,\, \tilde{s}\in [s_0, s_1-3],\cr
\|S_{\theta_{n}}\bar{e}_{\pm,n-1}\|_{s,T}\le
C\theta_n^{(s-\tilde{s})_+}\|\bar{e}_{\pm,n-1}\|_{\tilde {s},T}\le
C\theta_n^{(s-\tilde{s})_++L^2(\tilde{s})}\delta^2\Delta_n \qquad
\mbox{for}\,\, \tilde{s}\in [s_0, s_1-3]. \end{cases}
\end{equation}
Estimate $(\ref{f8.54})_2$ can be represented as the following three
cases.

{\it Case 1: $s=s_0$}. For this case, if we choose
$\tilde{s}=\alpha-3$, then $L^2(\tilde{s})=s_0-\alpha-1$, which
follows from  $(\ref{f8.54})_2$  that
$\|S_{\theta_{n}}\bar{e}_{\pm,n-1}\|_{s_0,T}\le C\delta^2
\theta_n^{s_0-\alpha-1}\Delta_n$;

{\it Case 2: $s=s_0+1$}. For this case, if we choose
$\tilde{s}=\alpha-2$, then $L^2(\tilde{s})=s_0-\alpha$, which
implies
$\|S_{\theta_{n}}\bar{e}_{\pm,n-1}\|_{s_0+1,T}\le C\delta^2
\theta_n^{s_0-\alpha}\Delta_n$;

{\it Case 3: $s\ge s_0+2$}. For this case, if we choose
$\tilde{s}=\alpha-4\ge s_0$, then $L^2(\tilde{s})=s_0-2-\alpha$,
which implies
$\|S_{\theta_{n}}\bar{e}_{\pm,n-1}\|_{s,T}\le
C\delta^2\theta_n^{s-\alpha-1}\Delta_n$.

In summary, we obtain
\begin{equation}\label{f8.58}
\|S_{\theta_{n}}\bar{e}_{\pm,n-1}\|_{s,T}\le
C\delta^2\theta_n^{s-\alpha-1}\Delta_n \qquad\mbox{for all}\,\,s\ge
s_0.
\end{equation}

Letting $\tilde{s}=\alpha$ in $(\ref{f8.54})_1$, we find
\begin{equation}\label{f8.59} \|h^\pm_n\|_{s,T}\le
C\delta^2\theta_n^{s-\alpha-1}\Delta_n \qquad\mbox{for all}\,\, s\ge
s_0,
\end{equation}
by using (\ref{f8.58}) in $(\ref{f8.41})_1$.

\vspace{.25in} \noindent {\bf Acknowledgements:} The research of
Gui-Qiang Chen was supported in part by the National Science
Foundation under Grants DMS-0505473, DMS-0426172, DMS-0244473, and
an Alexandre von Humboldt Foundation Fellowship. The research of
Ya-Guang Wang was supported in part  by a key project from the NSFC
under the Grant 10531020, a joint project from the NSAF and a
Post-Qimingxing Fund from the Shanghai Science and Technology
Committee under the Grant 03QMH1407. The second author would like to
express his gratitude to the Department of Mathematics of
Northwestern University (USA) for the hospitality, where this work
was initiated when he visited there during the Spring Quarter 2005.

\end{document}